%% file: announcements1.tex
\documentclass[11pt, a4paper]{amsart}
\usepackage[latin1]{inputenc}
\usepackage[T1]{fontenc}
\usepackage{amsmath, amssymb, amsthm}            
\usepackage{amstext, amsfonts, a4}

\usepackage[colorlinks=true]{hyperref} 
\hypersetup{urlcolor=blue, citecolor=red}
\usepackage[capitalise]{cleveref} 

\usepackage{color,transparent}
\usepackage{graphics,graphicx,subfigure}

\usepackage{tikz}
\usetikzlibrary{positioning,arrows,shapes,shadows}

\usepackage{comment}
\usepackage{float} 
\makeatletter

\def\myMRbibitem{\@ifnextchar[\my@lbibitem\my@bibitem}

\def\mybiblabel#1#2{\@biblabel{{\hyperref{http://www.ams.org/mathscinet-getitem?mr=#1}{}{}{#2}}}}

\def\myhyperanchor#1{\Hy@raisedlink{\hyper@anchorstart{cite.#1}\hyper@anchorend}}

\def\my@lbibitem[#1]#2#3#4\par{%
    \item[\mybiblabel{#2}{#1}\myhyperanchor{#3}\hfill]#4%
    \@ifundefined{ifbackrefparscan}{}{\BR@backref{#3}}%
    \if@filesw{\let\protect\noexpand\immediate
       \write\@auxout{\string\bibcite{#3}{#1}}}\fi\ignorespaces%
}

\def\my@bibitem#1#2#3\par{%
    \refstepcounter\@listctr
    \item[\mybiblabel{#1}{\the\value\@listctr}\myhyperanchor{#2}\hfill]#3%
    \@ifundefined{ifbackrefparscan}{}{\BR@backref{#2}}%
    \if@filesw\immediate\write\@auxout
        {\string\bibcite{#2}{\the\value\@listctr}}\fi\ignorespaces%
}

\makeatother

\newtheoremstyle{note}
  {5pt}
  {5pt}
  {\small}
  {}
  {\bfseries}
  {.}
  {.5em}
  {}

\theoremstyle{plain}
	\newtheorem{theorem}{Theorem} 
	\newtheorem{proposition}{Proposition}[section]        
	\newtheorem{problem}{Problem}[section]	    
	\newtheorem{corollary}{Corollary}[theorem]
	\newtheorem{conjecture}{Conjecture}[proposition]
	\newtheorem{claim}{Claim}[section]

        \newtheorem{scholium}[proposition]{Scholium}
        \newtheorem*{fact}{Fact}

\theoremstyle{definition}
	\newtheorem{definition}{Definition}[proposition]
        \newtheorem{example}{Example}[proposition]

\theoremstyle{remark}
	\newtheorem{remark}[proposition]{Remark}

	\newtheorem{question}[proposition]{Question}
        \newtheorem*{ack}{Acknowledgements}    

\numberwithin{equation}{section}         

\crefname{subsection}{\S}{\S\S} 
\Crefname{subsection}{\S}{\S\S}

\def\ogg~{{\rm \og}}   

             \def\NN{{\mathbb N}}    \def\RR{{\mathbb R}}        \def\ZZ{{\mathbb Z}}

  \def\cG{{\mathcal G}} \def\cM{{\mathcal M}}       \def\cC{{\mathcal C}}    \def\cU{{\mathcal U}} \def\cD{{\mathcal D}}    \def\cV{{\mathcal V}}         \def\cR{{\mathcal R}}    

\DeclareMathOperator{\Diff}{{Diff}}

\DeclareMathOperator{\ind}{ind}

\DeclareMathOperator{\dist}{dist}

\DeclareMathOperator*{\supp}{supp}

\def\vdashv{\vdash\!\dashv}

\def\Diff{\operatorname{Diff}}

\def\dist{\operatorname{dist}}
\def\Orb{\operatorname{Orb}}

\def\Het{\operatorname{Het}}
\def\Hom{\operatorname{Hom}}

\begin{document}
\title[$C^r$-generic dynamics close to homoclinic/heteroclinic points]{Steps towards a classification of $C^r$-generic dynamics close to homoclinic/heteroclinic points}
\author[N. Gourmelon]{N. Gourmelon}
\address{IMB. UMR 5251. \rm Universit\'e Bordeaux 1, 33405 Talence, France}
\email{Nicolas.Gourmelon "at" math.u-bordeaux1.fr}
\thanks{Supported by the IMB, Université de Bordeaux.}
\thanks{Partially supported by the ANR project \emph{DynNonHyp} BLAN08-2 313375}
\thanks{Partially supported by EU Marie-Curie IRSES "Brazilian-European partnership in Dynamical 
Systems" (FP7-PEOPLE-2012-IRSES 318999 BREUDS)}
\thanks{Partially supported by IMPA - Rio de Janeiro and by the Laboratoire de Mathématiques (UMR 8628), Faculté des Sciences d'Orsay, Université Paris-Sud.}

\keywords{Discrete-time systems; semilinear systems; bilinear systems; universal regular control}
\subjclass[2010]{37Cxx, 37Dxx, 37Gxx; 28Dxx, 34C23, 34C37, 34K18.}

\date\today

\begin{abstract}
We present here the first part of a program for a classification of the generic dynamics that may appear close to homoclinic and heteroclinic points, in the $C^r$ topologies, $r\geq 1$. This paper contains only announcements, formalism, and a few sketches of proofs. A forthcoming series of papers will present the  proofs in details.

The two prototypical examples of non-hyperbolic dynamics are homoclinic tangencies and heterodimensional cycles. Palis conjectured that they actually characterize densely non-hyperbolic dynamics. It is therefore important to understand what happens close to those bifurcations. 

A classical result of Newhouse~\cite{N2} from the seventies asserts that close to diffeomorphisms of a surface with a homoclinic tangency there is abundance of diffeomorphisms exhibiting infinitely many sinks or sources. This was extended to higher dimensions by Palis and Viana~\cite{PV}, in the nineties, assuming however that the tangency is sectionally dissipative, which implies in particular that it has index $1$. 
We propose a full generalization of the results of Newhouse, Palis and Viana, for both tangencies and cycles: close to a homoclinic tangency or to a heterodimensional cycle there is abundance of diffeomorphisms exhibiting infinitely many sinks or sources if and only if the dynamics is not volume-hyperbolic. This implies in particular a conjecture of Turaev for homoclinic tangencies. 

The dichotomy between volume-hyperbolicity and Newhouse phenomena actually holds in a much more general setting. It holds on "heteroclinic graphs", that is, on heteroclinic combinations of hyperbolic sets, under a generic combinatorial condition. 

An important result by Bonatti, Diaz, Pujals~\cite{BDP} states that if a homoclinic class is $C^1$-robustly without dominated splitting, then  nearby diffeomorphisms exhibit $C^1$-generically infinitely many sinks or sources. 
We obtain that this holds in $C^r$-regularity, for $r>1$, under the further assumption that non-dominations are obtained through so-called "mechanisms". This includes all the examples of robustly non-dominated homoclinic classes one can build with the tools known up to now.

We actually have a $C^r$-equivalent of a recent $C^1$-result of Bochi and Bonatti~\cite{BB} as we describe precisely the Lyapunov exponents along periodic points that may appear close to a homoclinic tangency or to a homoclinic class. 

The results of Newhouse, Palis and Viana were proven for the $C^r$ topologies, $r\geq 2$. Our results hold also in the $C^{1+\alpha}$ topologies.

\end{abstract}

\maketitle


\section{Introduction}
One wants to describe typical dynamical patterns, that is, the largest possible classes of dynamical systems. By large sets, we mean open and dense sets, or {\em residual} sets, that is, a countable intersection of open and sense subsets. In this paper we study the dynamics of $C^r$-generic diffeomorphisms, that is, the dynamics of residual subsets of diffeomorphisms of $\Diff^r(M)$, where $M$ is a compact Riemannian manifold and $r\geq 1$.

In the last 40 years, huge progress has been made in the classification of $C^1$-generic diffeomorphisms, building on a number of $C^1$-specific perturbative tools such as the closing lemma of Pugh, the ergodic closing lemma of Ma\~n\'e, the connecting lemma of Hayashi and the Franks' lemma, as well as a number of remarkable extensions (see~\cite{C1} for a survey). 

However, a residual set of $C^1$-diffeomorphisms may not intersect the set of $C^r$ diffeomorphism, for $r>1$. Therefore the dynamics of typical $C^1$-diffeomorphisms does not a priori indicate anything on the dynamics of typical $C^2$-diffeomorphisms. And indeed, Newhouse~\cite{N1} showed the existence of $C^2$-open sets of diffeomorphisms on surfaces presenting some dynamical phenomena (tangencies), while a result of Moreira~\cite{Mo} implies that these phenomena only appear for diffeomorphisms in a $C^1$-meager set.

In fact, the tools enumerated above are lacking in higher regularity: the $C^r$-closing and -connecting lemmas are still conjectures and the Franks' lemma does not hold in $C^r$-topology, for $r>1$. We present here a program whose aim is to provide a setting in which that shortage of tools is circumvented and a number of landmark results of $C^1$-generic dynamics are extended to the $C^r$-topology.
\medskip

Most of the results of this program are transposable to the volume preserving setting. 
\medskip

\subsection{Context} There exists an open set of diffeomorphisms -- now called hyperbolic diffeomorphisms -- whose dynamics has been extensively studied and are now well understood, thanks in particular to the spectral decomposition theorem of Smale (see for instance~\cite{BDV}).

Unfortunately, hyperbolic diffeomorphisms fail to be dense in $\Diff^1(M)$ when $\dim M\geq 3$~\cite{AS}, and in $\Diff^2(M)$ on surfaces~\cite{N1}. The two prototypical examples of non-hyperbolic dynamics are homoclinic tangencies and heterodimensional cycles.
Palis conjectured that they densely characterize open sets of non-hyperbolic dynamics:

\begin{conjecture}[Palis' density conjecture]
For $r\geq 1$, any diffeomorphism in $\Diff^r(M)$ is either approximated by hyperbolic ones or by diffeomorphisms exhibiting a homoclinic tangency or a heterodimensional cycle.
\end{conjecture}

One main goal of differentiable dynamical systems is therefore to describe what phenomena appear close to tangencies and cycles with the study of
\begin{itemize}
\item unfolding of homoclinic tangencies, H\'enon-like families and related problems of renormalization theory, which constitute whole areas of dynamical systems. 

\item unfolding of heterodimensional cycles: there are a number of advances in particular by Diaz, Rocha, Bonatti, Kiriki, in the $C^1$-topology, and in smooth topologies by Turaev, Gonchenko, Shil'nikov, among other authors.

%
\end{itemize}

Much of the existing literature is about what happens close to these two well-delimited dynamical mechanisms that yield non-hyperbolic dynamics, and such will it remain in the $C^r$ topologies, for $r>1$, as long as the absence of closing and connecting lemmas gives no hope of a global description of generic dynamics.

One aspect of our work is to propose a  setting that unifies all the known dynamical mechanisms that yield robustly non-hyperbolic dynamics,  that puts aside the problems posed by these lacking tools, and in which the practical dynamical consequences of the Franks' lemma will hold.

\subsection{Informal presentation of results}
 The definitions and the precise statements are given in~\cref{s.defsandstatements}.

\subsubsection{Newhouse phenomena close to homoclinic tangencies and heteroclinic cycles}
In celebrated papers of the seventies, Newhouse~\cite{N1,N2} showed that homoclinic tangencies on surfaces lead to  particularly wild dynamics: 
nearby any smooth diffeomorphism on a surface with a homoclinic tangency, there is a residual subset of an open set of diffeomorphisms, such that any element in it displays infinitely many sinks or sources. That open set is often called a {\em Newhouse region}. 

That result was extended in the nineties to higher dimension by Palis and Viana~\cite{PV} for saddles that satisfy a technical condition "sectional dissipativity" that implies in particular that the saddle has stable or unstable index $1$. Gonchenko, Shil'nikov and Turaev also explored generic dynamics close to certain heteroclinic cycles~\cite{GST1,GST3}. 

A natural obstruction for the appearance of sinks or sources close to a set is a weak form of hyperbolicity called "volume-hyperbolicity". Turaev conjectured a sharp dichotomy, which we state informally:

\begin{conjecture}[Turaev]\label{c.turaev}
Let $r\geq 2$. There is a Newhouse region $C^r$-close to a homoclinic tangency if and only if that tangency is not volume-hyperbolic.
\end{conjecture}

Our program leads to a proof of the conjecture and actually to a full generalization of the results of Newhouse and Palis-Viana for both homoclinic tangencies and heteroclinic cycles.  Let us state it informally:
\medskip

\noindent {\bf Theorem.} {\em Let $r>1$. There is a sharp dichotomy: a homoclinic tangency or a heteroclinic cycle is 
\begin{itemize}
\item either volume-hyperbolic,
\item or there is a Newhouse region $C^r$-close to it.
\end{itemize} }

\medskip

The theorem is precisely stated in~\cref{s.turaevconj}. While the theorems of Newhouse, Palis-Viana hold in $C^r$-topologies, $r\geq 2$, our results are valid in the $C^{1+\alpha}$ topologies thanks to a theorem by Crovisier and the author, as explained in \cref{s.C1alpha}.

In \cref{s.moregenerallambda}, we state a similar dichotomy in the general setting "heteroclinic graphs", that is, on heteroclinic and homoclinic combinations of hyperbolic sets, under a generic combinatorial condition on the graph. 

\subsubsection{Mechanisms}
A {\em basic set} for a difeomorphism $f$ is a compact  hyperbolic set, that is transitive and locally maximal. Its {\em index} (also {\em stable index}) is the dimension of its stable manifold. Bonatti~\cite{B} proposed to call "mechanism" any of the two following dynamical objects: 
\begin{itemize}
\item robust tangencies: a basic set  $K$ such that for any small $C^r$-perturbation of the dynamics, the stable and unstable laminations of the continuation of $K$ meet non-transversally,
\item robust cycles: a pair of basic sets $K$ and $L$ such that, for any small $C^r$-perturbation, the unstable lamination of one meets the stable lamination of the other.
\end{itemize}
He conjectured that there is an open and dense subset of diffeomorphisms that are either hyperbolic, or exhibit one of these two mechanism. 
These two are "mechanisms" in that they are well-delimited dynamical objects that involve only a finite number of basic sets and their respective local stable and unstable laminations. 
We propose to extend the list of mechanisms to include:
\begin{itemize}
\item periodic points with {\em non-simple Lyapunov spectrum}: periodic points at which the derivative of the first return map has a pair of eigenvalues (counted with multiplicity) with same modulus,
\item tangencies between different basic sets: a pair of basic sets of same index such that the stable lamination of one meets non-transversally the unstable lamination of the other. More generally: "bundle tangencies", where a plane of a strong-stable or center-stable distribution on the stable manifold of one basic set intersects in a non-generic position a plane of a strong-unstable or center-unstable distribution on the unstable manifold of the other basic set. 
\end{itemize}

\begin{figure}
\begin{center}
\def\svgwidth{130pt}           
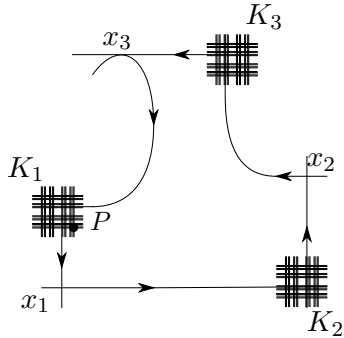
\end{center}
\caption{\small Cycle $(K_i,x_i)_{i\in \ZZ/3\ZZ}$ of basic hyperbolic sets containing the saddle $P$ and a tangency at $x_3$.}
\label{f.cyclemech}
\end{figure}

An invariant set is {\em dominated of index $i$}, if it admits a two-bundles dominated splitting where the first bundle has dimension $i$ (see~\cref{s.preldef}). 
We will say that the chain-recurrence class $\cC(P,f)$ of a saddle $P$ has {\em mechanically no domination of index $i$} if $P$ is in a cycle of basic hyperbolic sets that contains a mechanism -- either a periodic point with non-simple Lyapunov spectrum or a tangency/bundle tangency -- preventing domination of index $i$, as depicted in~\cref{f.cyclemech}. 
Here are the important facts underlying our definition of mechanisms:
\medskip

\noindent {\em 
\begin{enumerate}
\item A cycle of saddle points has no domination of index $i$ if and only if it has mechanically no domination of index $i$.
\item
All known examples of homoclinic classes/chain-recurrence classes robustly without domination of some index $i$, are obtained through mechanisms.
\end{enumerate}}

\medskip
Similarly to the conjecture of Bonatti, we conjecture that for a residual set of diffeomorphisms, non-dominations on the chain-recurrence class of any saddle $P$ are all mechanical. 
 This is presented with more details in \cref{s.mecha}; bundle tangencies are discussed in~\cref{s.bundtang}.

\subsubsection{Extensions of classical $C^1$-dynamics results to the $C^r$-topology through mechanisms}

By a famous result by Bonatti, Diaz and Pujals~\cite{BDP}, $C^1$-generically, if a homoclinic class has no domination, then it is the Hausdorff limit of a sequence of sinks or sources. 
We show that their result holds $C^r$-generically, for $r\geq 1$, under the further assumption that non-dominations are mechanical:
\medskip

\noindent {\bf Theorem}{\bf .} {\em Let $r\geq 1$. Given a $C^r$-generic diffeomorphism, if the chain-recurrence class of a saddle $P$ has  mechanically no domination of any index, then the homoclinic class of $P$ is the Hausdorff limit of a sequence of sinks or sources.}
\medskip

We refer the reader to \cref{s.Newhphdom} for definitions and formal statement.
As noted above, mechanisms include all the dynamical objects known to lead to robust non-domination, therefore  we get that:
\medskip 

\noindent {\em The Bonatti-Diaz-Pujals theorem holds for $r\geq 1$, for all the examples of homoclinic classes $C^r$-robustly without dominated splitting we can build with the tools we know of, to date.}
\medskip

Wen~\cite{W} and the author~\cite{G2} proved that a homoclinic class of a saddle $P$ has no domination of same index as $P$  if and only if there is a $C^1$-perturbation that creates a homoclinic tangency related to $P$. We prove the following:
\medskip

\noindent {\bf Theorem}{\bf .} {\em Let $r\geq 1$. If the chain-recurrence class of a saddle $P$ has mechanically no domination of same index as $P$, then there is a $C^r$-perturbation of the dynamics that creates a homoclinic tangency related to $P$.}

\medskip

Actually, we prove more generally that if the chain-recurrence class of a saddle $P$ has mechanically no-domination of index $i$, then one creates a bundle tangency of index $i$ related to $P$ by a $C^r$-perturbation. It is easy to deduce, with the techniques of~\cite{BCDG}, that if the homoclinic class of a saddle $P$ has  no domination of  index  $i$ then there is a $C^1$-perturbation that creates a bundle tangency of index $i$ related to $P$.
\medskip

Finally, the Bonatti-Diaz-Pujals results were generalized recently by Bochi and Bonatti~\cite{BB} who described the vectors of Lyapunov exponents that one gets along periodic orbits close to a homoclinic class, by $C^1$-perturbations, according to the indices for which there is no dominated splitting on the homoclinic class.  
In \cref{s.Lyapunov graphs}, we propose an extension to $C^r$-diffeomorphisms, using  mechanisms and their notion of Lyapunov graph. 

We deduce a full picture of the $C^r$-generic Lyapunov exponents along periodic orbits, close to homoclinic tangencies: those asymptotically realizable vectors of Lyapunov exponents (by $C^r$-perturbations) form a convex polytope, and they are locally generically simultaneously realized. More precisely, for a generic diffeomorphism in a nearby open set of dynamics, the vectors of Lyapunov exponents along periodic orbits are dense in a neighborhood of the polytope (see \cref{s.polytope}).

\subsection{About this paper}

All the results stated here are only announcements, although the main steps of the most difficult proofs are outlined. The theorems we announce are numbered more or less by increasing complexity of their statements. Such are the dependences:

\begin{center}
\begin{tikzpicture}
\tikzset{
doublearrow/.style={draw, thick, double distance=4pt, ->,>=angle 60},
thirdline/.style={draw, thin, ->, >=hooks}
}
\matrix [column sep=7mm, row sep=4mm] {    
\node (pr67) [ shape=rectangle,left=5pt,align=center] {\small \cref{p.gathering}};
   \node (pr68) [ shape=rectangle,right=5pt,align=center] {\small\cref{p.robustizing}};& &
  \node (th12) [shape=rectangle,left=5pt,align=center] {\small\cref{t.bobocr}};
  \node (pr7) [right=5pt,align=center] {\small\cref{p.geninfcon}};  \\
     \node (th11) [draw, very thick,shape=rectangle,align=center] {\small\cref{t.homcor}: robustization of\\ \small non-dominations on\\ \small heteroclinic graphs};&
   \node (th13) [shape=rectangle,align=center] {\small\cref{t.Lyaptang}};  &
   \node (cor12) [draw, very thick, shape=rectangle,align=center] {\small\cref{c.Lyapgraph}:  \\ \small a $C^r$-version of \cite{BB}\\ \small on heteroclinic classes}; \\
  \node (th2) [shape=rectangle,rectangle split, rectangle split parts=2,align=center,below left=5pt and -5pt] {$\begin{cases}\mbox{\small \cref{t.robustization}} \\ \mbox{\small \cref{t.robustizationalpha} } (C^{1+\alpha})\end{cases}$};  & 
         \node (th4) [shape=rectangle,align=center] {\small \cref{t.homoclnodommoregen}}; &

  \node (th3) [ shape=rectangle, left=5pt,align=center] {\small \cref{t.comprehensible}};  
    \node (th9) [ shape=rectangle,right=10pt,align=center] { \small\cref{t.mechnovolhyp}};\\
  };
\node (th10) [draw,very thick, shape=rectangle,align=center, below= 15 pt of th4] {\small\cref{t.volhyp}: Newhouse ph./Volume Hyperbolicity \\ dichotomy on heteroclinic graphs};
  
 \node (th1) [shape=rectangle,rectangle split, rectangle split parts=2,align=center, below left = 15pt and -100pt of th10] {$\begin{cases}\mbox{\small\cref{t.lambda}} \\ \mbox{\small\cref{t.lambdaalpha} } (C^{1+\alpha})\end{cases}$};
  \node (cor10) [shape=rectangle,align=center, below right = 15pt and -100pt of th10] {\small\cref{c.homcor}};


\draw[double,double distance=2pt,-implies] (th12) -- ([xshift=-35pt]cor12.north);
\draw[double,double distance=2pt,-implies] (pr7) -- ([xshift=+42.5pt]cor12.north);
\draw[double,double distance=2pt,-implies] (cor12) -- (th13);

\draw[double,double distance=3pt,-implies] ([xshift=-33pt]cor12.south) -- (th3);
\draw[-implies] ([xshift=-33pt]cor12.south) -- (th3);

\draw[double,double distance=3pt,-implies] ([xshift=-51pt]th11.south) -- (th2);
\draw[-implies] ([xshift=-51pt]th11.south) -- (th2);

\draw[double,double distance=3pt,-implies] ([xshift=+37pt]cor12.south) -- (th9);
\draw[-implies] ([xshift=+37pt]cor12.south) -- (th9);

\draw[double,double distance=3pt,-implies] (th9) -- (th3);
\draw[-implies] (th9) -- (th3);

\draw[double,double distance=2pt,-implies] (th3) -- (th4);
\draw[double,double distance=2pt,-implies] (th11) -- (th13);
\draw[double,double distance=2pt,-implies] (th11.south east) -- (th4.north west);
\draw[double,double distance=2pt,-implies] (pr67) -- ([xshift=-42pt]th11.north);
\draw[double,double distance=2pt,-implies] (pr68) -- ([xshift=+42pt]th11.north);
\draw[double,double distance=2pt,-implies] (th9) -- (th10.north east);
\draw[double,double distance=2pt,-implies] ([xshift=20pt]th11.south) -- ([xshift=-70pt]th10.north);

\draw[double,double distance=3pt,-implies] ([xshift=-81pt]th10.south) -- ([xshift=-10pt]th1.north);
\draw[-implies] ([xshift=-81pt]th10.south) -- ([xshift=-10pt]th1.north);

\draw[double,double distance=3pt,-implies] ([xshift=+60pt]th10.south) -- ([xshift=10.5pt]cor10.north);
\draw[-implies] ([xshift=+60pt]th10.south) -- ([xshift=10.5pt]cor10.north);
\end{tikzpicture}

{\small
Legend: \begin{tabular}{cl}
 $A \Rrightarrow B$ & : $B$ is a direct consequence of $A$ \\
 $A \Rightarrow C \Leftarrow B$ & : ($A$ and $B$) imply $C$.
 \end{tabular}}
 \end{center}

\section{Definitions and statement of results}\label{s.defsandstatements}
We now give formal definitions and state precisely the simple versions of our results. The introduction of the new notion of bundle tangency is postponed to \cref{s.mechnodom} since the larger part of these results can be understood without it. The $C^{1+\alpha}$ theorems are stated in \cref{s.C1alpha}.
 The most general theorems, involving heteroclinic combinations of basic sets, need an adequate formalism and are stated in \cref{s.moregenerallambda}. Finally, the results about the description of Lyapunov exponents along periodic orbits are postponed to \cref{s.Lyapunov graphs}.

\subsection{preliminary definitions}\label{s.preldef}

As we noted in the introduction, since a  classification of the global dynamics of $C^r$-generic diffeomorphisms seems out of reach, for $r>1$, most of the existing literature considers an invariant non-hyperbolic compact invariant set $\Lambda$ with a well-defined dynamical mechanisms (almost always a homoclinic tangency or a heterodimensional cycle), and studies what happens in a neighborhood of $\Lambda$ for generic perturbations of the dynamics.

We want to explore this question systematically for a widest possible class of sets $\Lambda$.
For simplicity, we propose the following:

\begin{definition}
Given a diffeomorphism $f\in \Diff^r(M)$ and an $f$-invariant set $\Lambda$, we will say that {\em Newhouse phenomena occur $C^r$-close to $\Lambda$} if the following holds:  for any open set $U$ containing $\Lambda$, there exists an open set $\cU\subset \Diff^r(M)$ such that 
\begin{itemize}
\item  $f$ is in the closure of $\cU$,
\item $C^r$-generic diffeomorphisms in $\cU$ have infinitely many sinks or sources, by restriction to $U$.
\end{itemize}
\end{definition}

The problem we deal with in the first part of our program may be formulated as follows:
\begin{problem}\label{p.newh}
Given an $f$-invariant set $\Lambda$, characterize the occurrence of Newhouse phenomena $C^r$-close to $\Lambda$ by the dynamical properties of the restriction $f_{|\Lambda}$.
\end{problem}

One may of course want to characterize all kinds of dynamical phenomena:

\begin{problem}\label{p.general}
Given an $f$-invariant set $\Lambda$, characterize the dynamical phenomena (for instance homoclinic tangencies, wild dynamics, universal dynamics, viral dynamics, etc...) that occur $C^r$-close to $\Lambda$ by the dynamical properties of the restriction $f_{|\Lambda}$.
\end{problem}

Let us now define a few of the "characteristic" dynamical properties of $f_{|\Lambda}$.

An invariant set $\Lambda$ is {\em hyperbolic} if there is a $Df$-invariant continuous splitting $TM_{|\Lambda}=E^s\oplus E^u$ of the tangent bundle such that there is some positive iterate of the dynamics that contracts (resp. dilates) by a factor $\lambda>1$ all nonzero vector of $E^s$ (resp. $E^u$). The {\em index} of $\Lambda$ is the dimension of its stable bundle.

The notion of hyperbolic splitting is naturally weakened to that of dominated splitting: a $Df$-invariant continuous splitting $TM_{|\Lambda}=E_1\oplus \ldots \oplus E_k$ into constant-dimensional vector bundles is {\em dominated} if there exists $\lambda>1$ and an iterate $N>0$ such that for any point $x\in \Lambda$, for any pair of unit vectors $u,v\in E_i, E_{i+1}$ in consecutive bundles and above $x$, we have
$$\|Df^N(u)\|<\lambda^{-1}\|Df^N(v)\|.$$
A {\em domination of index $i$} on an $f$-invariant set $\Lambda$ is a dominated splitting $E\oplus F$ such that the bundle $E$ has dimension $i$. 

An invariant set $\Lambda$ of a diffeomorphism $f$ is {\em volume-hyperbolic} if there exists a dominated splitting $T_\Lambda M= E^{vs}\oplus E^{c}\oplus E^{vu}$
such that the bundles $E^{vs},E^{vu}$ are non-trivial and such that volume is contracted on $E^{vs}$ and expanded on $E^{vu}$ by a factor $\lambda>1$ by some iterate of the dynamics, that is, for some  $N$, for all $x\in \Lambda$, the determinants $\det(Df^N_{|E^{vs}_x})$ and $\det(Df^{-N}_{|E^{vu}_x})$ are less than $\lambda^{-1}$.

\begin{remark}
If a set $\Lambda$ is volume-hyperbolic, clearly no sink nor source will appear by $C^1$ perturbations in a neighborhood of $\Lambda$.
\end{remark}

\begin{remark}
By~\cite{G1}, there are adapted metrics to dominated splittings: one may change the Riemannian metric on $M$ to an equivalent Riemannian metric and assume that $N=1$ in the definition of domination.
\end{remark}

\subsection{\cref{p.newh} for homoclinic tangencies and heterodimensional cycles}\label{s.turaevconj}
A {\em homoclinic tangency} is a saddle point $P$ whose stable and unstable manifold meet non-transversely at a point $x$. We say that we have a homoclinic tangency {\em related to $P$}. By abuse of terminology, we call {\em homoclinic tangency} the compact set 
\begin{align*}\Lambda&=\overline{\Orb(x)}\\
&=\Orb(P)\cup \Orb(x).
\end{align*}
 A {\em heteroclinic cycle} is a pair of saddle points $P,Q$ connected through two heteroclinic points $x\in W^s(P)\cap W^u(Q)$ and  $y\in W^u(P)\cap W^s(Q)$. Similarly we call {\em heteroclinic cycle} the set 
\begin{align*}
\Lambda&=\overline{\Orb(x)\cup \Orb(y)}\\
	&= \Orb(P)\cup \Orb(Q) \cup \Orb(x)\cup \Orb(y).
\end{align*}
A {\em heterodimensional cycle} is a particular case of heteroclinic cycle where $P$ and $Q$ have different indices.
Let us restate the classical results  of Newhouse in dimension $2$:
\medskip

\noindent {\bf Theorem} (Newhouse~\cite{N1,N2}). {\em Let $r\geq 2$ and let $\Lambda$ be a homoclinic tangency on a surface. Then Newhouse phenomena occur $C^r$-close to $\Lambda$.}
\medskip

In a remarkable paper, Palis and Viana extended this to higher dimension under a technical hypothesis: the $f$-invariant set $\Lambda$ is called {\em sectionally dissipative} if, either for $f$ or its inverse, volume on any $2$-plane $P$ in $T_\Lambda M$ is exponentially contracted by positive iterations. In particular, the saddle must have index $1$ or coindex $1$. 
 \medskip

 \noindent {\bf Theorem} (Palis-Viana~\cite{PV}). {\em Let $r\geq 2$ and let $\Lambda$ be a sectionally dissipative homoclinic tangency. Then Newhouse phenomena occur $C^r$-close to $\Lambda$. }
\medskip
 
Dealing with higher index saddles adds difficulty. As seen in the introduction, Turaev conjectured that volume-hyperbolicity is the only obstruction to the occurrence of Newhouse phenomena $C^r$-close to a homoclinic tangency.
Some particular cases of heteroclinic cycles were extensively studied by Shil'nikov, Gonchenko and Turaev~\cite{GST1,GST3}. 

We now restate formally our main theorem, which fully generalizes the Newhouse, Palis-Viana theorems. It implies  \cref{c.turaev} and solves \cref{p.newh} for both homoclinic tangencies and heteroclinic cycles:
  
\begin{theorem}\label{t.lambda}
Let $r\geq 2$ be an integer. Let $\Lambda$ be a homoclinic tangency or a heteroclinic cycle for a diffeomorphism $f\in \Diff^r(M)$. One has a sharp dichotomy:
\begin{itemize}
\item either $\Lambda$ is volume-hyperbolic,
\item or  Newhouse phenomena occur $C^r$-close to $\Lambda$
\end{itemize}
\end{theorem}

The first step to prove \cref{t.lambda} is to turn robust, that is, resistant to small perturbations, the lack of domination on $\Lambda$. 
Recall first that the {\em homoclinic class} of a saddle $P$ is defined by the closure of the transverse intersections of the stable and unstable manifolds of $P$:
$$\Hom(P,f)=\overline{W^s(f,P)\pitchfork W^u(f,P)}.$$
Here is the robustization theorem:

\begin{theorem}\label{t.robustization}
Let $r\geq 2$ be an integer.  Let $\Lambda$ be a homoclinic tangency or a heteroclinic cycle containing a saddle $P$  for a diffeomorphism $f\in \Diff^r(M)$. Assume that $\Lambda$ has no dominated splitting of any index in a set $I\subset \NN$.

Then one can make simultaneously robust  those non-dominations by a $C^r$-perturbation: in any neighborhood $\cV\subset \Diff^r(M)$ of $f$, there exists an open set $\cU\subset \cV$ of diffeomorphisms $g$ such that the homoclinic class $\Hom(P_g,g)$ of the continuation $P_g$ of $P$ has no domination of any index $i\in I$.
\end{theorem}

The theorem that we actually prove is  more precise in its conclusions: the robust non-dominations are obtained through the mechanisms evoked in the introduction and defined in \cref{s.mecha}. 
We will precisely use these mechanisms to eventually create sinks or sources.

\begin{figure}[H]\label{f.grapheeulerienselles}
\begin{center}
\def\svgwidth{130pt}           
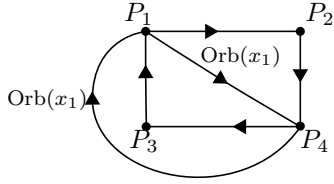
\end{center}
\caption{\small Example of a set $\Lambda$, "Eulerian" union of saddles and heteroclinic orbits. There are Newhouse phenomena $C^r$-close to this set if and only if it is not volume hyperbolic.}
\label{f.grapheeul}
\end{figure}

In fact, both theorems hold for more general sets $\Lambda$:

\medskip

\noindent {\bf Addendum to \cref{t.lambda,t.robustization}:}
{\em both theorems still hold if the homoclinic tangency or heterodimensional cycle $\Lambda$ is replaced by a union 
$$\Lambda=\bigcup_{1\leq j\leq n}\Orb(P_j) \cup \bigcup_{1\leq j\leq m}\Orb(x_j)$$
where each $P_j$ is a saddle point, so that:
\begin{itemize}
\item the $\alpha$-limit (resp. $\omega$-limit) of each point $x_j$ is one of the orbits $\Orb(P_k)$,
\item if we define an oriented graph $\Gamma$ so that each $\Orb(P_j)$ is a vertex and each $\Orb(x_j)$ is an edge from $\alpha(x_j)$ to $\omega(x_j)$, then that graph is Eulerian, that is, there is a closed path that goes through each edge once and only once (see \cref{f.grapheeulerienselles}).
\end{itemize}}

\begin{remark}
We point out that these theorems
\begin{itemize}
\item also hold in $C^{1+\alpha}$ topologies; this is discussed in \cref{s.C1alpha}.
\item hold for even more general sets $\Lambda$, namely, for heteroclinic combinations of basic sets under a similar Eulerian condition, which happens to be generic. This is discussed in \cref{s.moregenerallambda}.
\end{itemize}
\end{remark}

\begin{remark}
The Eulerian condition in the addendum is essential. It can not be replaced by a weaker strong connectivity condition, as the counter-example of \cref{s.counterex} shows (see also \cref{f.counterex3}).
\end{remark}

\subsection{\cref{p.newh}  for chain-recurrence classes}

The results stated here are $C^r$-extensions of theorems of $C^1$-generic dynamics~\cite{BDP,W,G2,BB}. They are needed to deduce \cref{t.lambda} from \cref{t.robustization}. 

\subsubsection{Chain-recurrence classes}
In dissipative dynamics, all the complexity of the dynamics (entropy, etc...) is visible by restriction to a set of points who present a certain form of recurrence. Conley~\cite{Co} proposed a weakest form of recurrence, now widely used: chain-recurrence, that is recurrence through pseudo-orbits. 

An {\em $\epsilon$-pseudo orbit} is a sequence $(x_i)_{1\leq i\leq n}$ in $M$ such that the distance from $f(x_i)$ to $x_{i+1}$ is less than $\epsilon$. 
We write $x${\SMALL$\vdashv$}$y$ if, for all $\epsilon>0$, there exist $\epsilon$-pseudo-orbits of length $>1$ from $x$ to $y$ and from $y$ to $x$. The {\em chain-recurrent set} is the set $\cC(f)$ of points $x$ such that $x${\SMALL$\vdashv$}$x$. The relation {\SMALL$\vdashv$} is an equivalence relation on $\cC(f)$. The {\em chain-recurrence class $\cC(x,f)$} of $x\in \cC(f)$ is the equivalence class of $x$.

 An $f$-invariant set $\Sigma\subset \cC(f)$ is called {\em chain-recurrent} if for all $x,y\in \Sigma$, for all $\epsilon>0$, there exist $\epsilon$-pseudo-orbits in $\Sigma$ from $x$ to $y$ and from $y$ to $x$. 
\medskip

If a point $P$ is a saddle then the chain-recurrence and homoclinic classes of $P$ will actually be identified to pointed sets:
\begin{align*}\cC(P,f)&:=\bigl(\cC(P,f),P\bigr)\\
\Hom(P,f)&:=\bigl(\Hom(P,f),P\bigr).
\end{align*}
The continuations of the chain-recurrence and homoclinic classes of a saddle $P$ for $f$ are therefore well-defined, for any diffeomorphism $g$ in a small neighborhood of $f$.

\subsubsection{Mechanisms}\label{s.mecha}
In~\cite{B}, Bonatti proposed to gather under the term "mechanisms" dynamical objets that are known to robustly prevent hyperbolicity. Similarly, we propose a simple and comprehensive list of "mechanical obstructions to domination": it gathers all the known ways to build homoclinic classes/chain-recurrence classes that are robustly non-dominated of some index $i$. 

We need first a few definitions. A {\em basic set} is a transitive locally maximal invariant hyperbolic set. A {\em cycle of basic sets} is a cyclic sequence $(K_\ell,x_\ell)_{\ell\in \ZZ/n\ZZ}$ where the $K_\ell$ are basic sets -- not necessarily of same index -- and $x_i$ is an intersection point of the unstable lamination of $K_\ell$ and the stable one of $K_{\ell+1}$.  We call {\em mechanical obstruction} to domination of index $i$ any of the following:
\begin{itemize}
\item a periodic point $Q$ of period $p$ such that the eigenvalues $\lambda_1,...,\lambda_d$ of the derivative $Df^p$ at $Q$, counted with multiplicity and ordered by increasing moduli, satisfy $|\lambda_i|=|\lambda_{i+1}|$.
\item a {\em tangency}, that is, a pair of basic sets $K,L$ of index $i$ whose stable and unstable bundles meet non-transversely, or more generally, a {\em bundle tangency of index $i$}, which we define precisely in \cref{s.bundtang}: 

a plane of a strong-unstable or center-unstable distribution on the unstable manifold of $K$ intersects non-generically a plane of a strong-stable or center-stable distribution on the stable manifold of $L$, so that domination of index $i$ is prevented.
\end{itemize}
We say that the chain-recurrent class $\cC(P,f)$ of a saddle $P$ has {\em mechanically no domination of index $i$} if there is a cycle of basic sets that contains $P$ and a mechanical obstruction to domination of index $i$. 

\begin{remark}
Two chain-recurrence classes $\cC(P,f)$ and $\cC(Q,f)$ may coincide as subsets of $M$, but one may have mechanically no domination of some index while the other not. This is not an abuse of terminology if, as said previously, we consider chain-recurrent classes to be pointed sets.
\end{remark}

\begin{fact}
In all known examples, $C^1$- or $C^r$-robust obstructions to domination on the chain-recurrent class $\cC(P,f)$ of a saddle (see~\cite{N1,N2,PV,GST1,GST3,BD}, among other works) are locally openly and densely mechanical.
\end{fact}

The same way Bonatti~\cite{B} conjectured that robust tangencies and heterodimensional cycles  characterize openly and densely non-hyperbolic dynamics, we conjecture that mechanisms encompass all the generic obstructions to dominated splittings:

\begin{conjecture}[non-dominations are generically mechanical]
There is a residual subset $\cR$ of $C^r$-diffeomorphisms $f$ such that if the chain-recurrent class $\cC(P,f)$ of a saddle $P$ is not dominated of index $i$, then it has $C^r$-robustly mechanically no domination of index $i$.
\end{conjecture}

This conjecture is partially proved in the $C^1$-topology as a consequence of the works of \cite{BC1} and \cite{BCDG}.

\subsubsection{Newhouse phenomena versus domination}
\label{s.Newhphdom}

To understand the theorems stated here in their fullest extent, one would need to read the lengthy definition of bundle tangency in \cref{s.mechnodom}. We postponed that definition since the larger part of these results can be understood without it. 

\medskip

 \cref{p.newh} for homoclinic classes and in the $C^1$-topology is already quite well understood, thanks to the existing perturbative tools: closing, ergodic closing, connecting, Franks' lemmas and their multiple extensions.
In a landmark paper of $C^1$-generic dynamics, Bonatti, Diaz and Pujals~\cite{BDP} showed that, $C^1$-generically, if a homoclinic class $\Lambda=\Hom(P,f)$ has no domination of any index, then it is the Hausdorff limit of a sequence of sinks and sources. 
We propose an extension to the $C^r$-topology, for $r\geq 1$; we show that the result holds for all the known examples of $C^r$-generic homoclinic classes, to date.

\begin{theorem}\label{t.comprehensible}
For any integer $r\geq 1$, for any diffeomorphism $f$ in a residual subset $\cR \subset \Diff^r(M)$, 
if the chain-recurrent class $\cC(P,f)$ of a saddle $P$ has mechanically no domination of any index, that is:

 if, for any index $1\leq i<d$, there is a cycle of basic sets that contains the saddle $P$ and either
\begin{itemize}
\item a periodic point $Q$ of period $p$ such that the eigenvalues $\lambda_1,...,\lambda_d$ of the derivative $Df^p$ at $Q$, counted with multiplicity and ordered by increasing moduli, satisfy $|\lambda_i|=|\lambda_{i+1}|$,
\item or a non-transverse intersection point of the stable and unstable laminations of to two basic sets $K,L$ of index $i$ (possibly $K=L$),
\item or, more generally a bundle tangency of index $i$, as defined later in~section~\ref{s.bundtang},
\end{itemize}
then the homoclinic class $\Hom(P,f)$ is the Hausdorff limit of a sequence of sinks or sources.
\end{theorem}

Examples of chain-recurrent classes having robustly mechanically no domination can be easily built in many ways, combining together hyperbolic sets that may contain saddles with complex eigenvalues, using blenders as in~\cite{BD}, or large Hausdorff dimension horseshoes (in the $C^r$-topologies) to create robust tangencies or robust heterodimensional cycles. 

With the fact stated in \cref{s.mecha}, \cref{t.comprehensible} implies the following:
\medskip

\noindent {\em In all known examples of chain-recurrent class without dominated splitting, there is $C^r$-generically a sequence of sinks or sources that Hausdorff converges to the corresponding homoclinic class.}
\medskip

\cref{t.comprehensible} implies that if a chain-recurrent class is robustly mechanically without domination of any index, then the homoclinic class is locally generically the Hausdorff limit of a sequence of sinks or sources. The hypotheses can be weakened; it is indeed enough to ask that mechanical non-dominations be locally dense, for each index $i$:

\begin{theorem}\label{t.homoclnodommoregen}
Let $r\geq 1$ be an integer and let $P$ be a saddle point for a diffeomorphism $f$. Let $\cU\subset \Diff^r(M)$ be an open neighborhood of $f$ such that 
\begin{itemize}
\item the continuation $P_g$ of $P$ is defined for $g\in \cU$,
\item for all $1\leq i< d$, there is a dense subset $\cD_i \subset \cU$ of diffeomorphisms $g$ such that the chain-recurrent class $\cC(P_g,g)$ has mechanically no domination of index $i$. 
\end{itemize}
Then there is a residual set $\cR\subset \cU$ of diffeomorphisms $g$ such that the homoclinic class $\Hom(P_g,g)$ is the Hausdorff limit of a sequence of sinks or sources.
\end{theorem}

We also get a $C^r$-version of the most general statement of~\cite{BDP}, which says that if a homoclinic class is $C^1$-robustly not volume-hyperbolic, then Newhouse phenomena occur $C^1$-close to it. We refer the reader to \cref{s.volhyp} for a precise statement.

\subsubsection{Tangencies versus domination}
Wen~\cite{W} and the author~\cite{G2} proved that if a homoclinic class of a saddle $P$ has no domination of same index as $P$ then there is a $C^1$-perturbation that creates a homoclinic tangency related to $P$, that is,  the stable and unstable manifolds of $P$ meet non-transversely after perturbation. We extend it to the $C^r$-topologies, again under the assumptions that non-dominations are mechanical:

\begin{theorem}\label{t.nodomtang}
 Let $r\geq 1$ be an integer and let $f\in\Diff^r(M)$. Let $P$ be a saddle orbit for $f$ such that its chain-recurrent class $\cC(P,f)$ has mechanically no domination of same index as $P$.
 
 In any $C^r$-neighborhood $\cU\subset \Diff^r(M)$ of $f$ there exists a diffeomorphism $g$ such that $P$ is again a saddle for $g$, with a homoclinic tangency related to it.
\end{theorem}

We actually prove the more general theorem:

\begin{theorem}\label{t.nodombundtang}
 Let $r\geq 1$ be an integer and let $f\in\Diff^r(M)$. Let $P$ be a saddle for $f$ such that its chain-recurrent class $\cC(P,f)$ has mechanically no domination of index $i$. 
 
 In any $C^r$-neighborhood $\cU\subset \Diff^r(M)$ of $f$ there exists a diffeomorphism $g$ such that $P$ is a saddle for $g$, with a bundle tangency of index $i$ related to it.
\end{theorem}

With the techniques of~\cite{BCDG}, one deduces that if the homoclinic class of a saddle $P$ has no domination of  index $i$ then there is a $C^1$-perturbation that creates a bundle tangency of index $i$ related to $P$.

\subsubsection{Description of Lyapunov exponents  for periodic points close to homoclinic classes}
The results of Bonatti-Diaz-Pujals~\cite{BDP} were generalized recently by Bochi and Bonatti~\cite{BB} who described the vectors of Lyapunov exponents that one gets along periodic orbits close to a homoclinic class, by $C^1$-perturbations, according to the indices for which there is no dominated splitting on the homoclinic class.  

We propose a $C^r$-version of~\cite{BB}: using their idea of Lyapunov graph, we describe the vectors of Lyapunov exponents that one gets along periodic orbits close to a generic homoclinic class, by $C^r$-perturbations, according to the indices for which there is mechanically no dominated splitting. 

With \cref{t.robustization}, we deduce a full picture of the $C^r$-generic Lyapunov exponents along periodic orbits, close to homoclinic tangencies: the vectors of Lyapunov exponents of the periodic points that may appear asymptotically close to a homoclinic tangency,  by $C^r$-perturbations, form a convex polytope. 

Moreover, for $C^r$-generic diffeomorphisms $g$ in a nearby open set $\cU\subset \Diff^r(M)$, the vectors of Lyapunov exponents for periodic points of $g$ form a dense set of a (small) neighborhood of that polytope. This is presented in \cref{s.Lyapunov graphs}.

\section{Bundle tangencies and heteroclinic classes}\label{s.mechnodom}
We give the definition of bundle tangency, needed to understand \cref{t.comprehensible,t.homoclnodommoregen,t.nodomtang,t.nodombundtang} to their full extent. We also discuss the notion of heteroclinic relation which is a natural extension of the notion of homoclinic relation, and which gives a clearer understanding of our work.

\subsection{Bundle tangencies}\label{s.bundtang}
We define a notion that generalizes homoclinic tangencies.
Note first that if an invariant set $K$ has a dominated splitting $TM_K=E\oplus F$ of some index $i$, that is, the bundle $\dim(E)$ has constant dimension $i$, then the center unstable bundle $F$ is uniquely extended to a  bundle $\tilde{F}$ above the unstable set 
$$W^u(K)=\{x\in M, \alpha(x)\subset K\}$$
 that is $Df$-invariant and locally continuous, that is, continuous by restriction to a local unstable set
 $$W^u_\epsilon(K)=\{x\in W^u(K) \mbox{ s.t. } \forall n\in \NN, \dist(f^{-n}(x), K)<\epsilon\},$$
 for some $\epsilon>0$.  Symmetrically, the center stable bundle $E$ is uniquely extended to a $Df$-invariant and locally continuous bundle $\tilde{E}$ above the stable set $W^s(K)$.

\begin{definition}
A {\em bundle tangency} is the data of
\begin{itemize}
\item a point $x\in M$
\item two dominated splittings $TM_{|\alpha(x)}=E^\alpha\oplus F^\alpha$ and $TM_{|\omega(x)}=E^\omega\oplus F^\omega$ on the $\alpha$- and $\omega$-limits of $x$, such that 
the extended bundles $\tilde{F}^\alpha$ and $\tilde{E}^\omega$ are not in generic position at $x$, that is, both  following conditions are met: 
\begin{align}
&d_T=\dim(\tilde{F}^\alpha_x\cap \tilde{E}^\omega_x)>0,\label{e.effezf1}\\
&T_xM\neq \tilde{F}^\alpha_x\oplus \tilde{E}^\omega_x \label{e.effezf2}.
\end{align} 
\end{itemize}
\end{definition}

We also say that we have a {\em bundle tangency at $x$, associated to the dominated splttings $E^\alpha\oplus F^\alpha$ and $E^\omega\oplus F^\omega$}.
We call $d_T$ the {\em dimension} of the bundle tangency (the dimension of a bundle tangency does {\em not} correspond to its degree of degeneracy).

Note that \cref{e.effezf1,e.effezf2} force that 
\begin{align}
i_\alpha+d_T&>i_\omega, \label{e.zfef3}
\end{align}
where $i_\alpha$ and $i_\omega$ are the respective indices of the dominated splittings $E^\alpha\oplus F^\alpha$ and $E^\omega\oplus F^\omega$.

\begin{scholium}\label{s.gentangency}
If we have a bundle tangency of dimension $d_T$ at $x$ associated to dominated splittings of indices $i_\alpha$ and $i_\omega$, then  there is no dominated splitting of any index $i$ such that 
$$i_\omega-d_T<i<i_\alpha+d_T$$
on the compact invariant set $\alpha(x)\sqcup \Orb(x) \sqcup \omega(x)$.
\end{scholium}

We call those integers $i$ the {\em indices} of the bundle tangency. By \cref{e.effezf1,e.zfef3} the set of indices is non-empty.  One deduces from the scholium that the point $x$ does not belong to the limit sets $\alpha(x)$ and $ \omega(x)$.

A bundle tangency {\em related to a saddle $P$}  is a bundle tangency where $\alpha(x)=\omega(x)=\Orb(P)$.

\begin{example}
A homoclinic tangency related to a saddle point $P$ is a bundle tangency related to $P$ and associated to the hyperbolic splitting $T_PM=E^s_{\Orb(P)}\oplus E^u_{\Orb(P)}$ along the orbit of $P$. Conversely, if the index of $P$ is also an index of a bundle tangency related to $P$, then that bundle tangency is a homoclinic tangency.
\end{example}

\begin{example}
If a saddle point $P$ of stable index $2$ in $\RR^3$ has eigenvalues 
\begin{align*}
0<\lambda_1<\lambda_2< 1< \lambda_3,
\end{align*} 
then there is a strong stable foliation of dimension $1$ on the stable manifold of $P$. If the unstable manifold meets the stable one tangently to that foliation at some point $x$, then there is at $x$ a bundle tangency of indices $1$ and $2$ related to $P$.
\end{example}

A {\em robust bundle tangency of index $i$} is a pair of basic sets $K$ and $L$ for a diffeomorphism $f$ such that for any neighboring diffeomorphism $g$, the stable and unstable manfifolds of the continuations $K_g$ and $L_g$ meet at a bundle tangency of index $i$. The first examples of robust homoclinic tangencies are due to~\cite{N1,PV,Ro}. But it needs not be a homoclinic tangency; we mention indeed the following example without expliciting the construction:

\begin{example}
Embedding for a instance Plykin attractor $K$ into a  normally hyperbolic manifold, it is easy to build a $C^1$-robust bundle tangency related to $K$ that is not a homoclinic tangency. This is a robust obstruction to domination on the homoclinic class of $K$ that is not a tangency nor a pair of complex conjugate eigenvalues.
\end{example}

Nevertheless, one can prove in this example that some $C^r$-perturbation will create a saddle orbit homoclinically related to $K$, along which there is a pair of complex conjugate eigenvalues of the corresponding index. 

\subsection{Heteroclinic classes}\label{s.hetclas}
The notion we define here is natural and facilitates the statement of our results. 
We say that two basic sets $K$ and $L$ are {\em heteroclinically related} if and only if they belong to some common cycle of basic sets. This is an equivalence relation on basic sets. Recall that the homoclinic class of a saddle point $P$ is the closure of the transverse intersections of the stable and unstable manifolds of $P$.

\begin{definition}
Given a diffeomorphism $f$, the {\em heteroclinic class} $\Het(K,f)$ of a basic set $K$ is the union of the basic sets $K_\ell$ and the points $x_\ell$, such that $K,K_\ell,x_\ell$ are part of a common cycle of basic sets.
\end{definition}

\begin{remark}
\begin{itemize}
\item A heteroclinic class is chain-recurrent.
\item The closure of a heteroclinic class contains the homoclinic class of each of its saddle points. 
\item In general, it is not the union of these homoclinic classes.
\end{itemize}
\end{remark}

\begin{proposition}\label{p.hethom}
 For any integer $r\geq 1$, there is a residual subset of $C^r$-diffeomorphisms $f$ such that 
\begin{enumerate}
\item for any saddle point $P$, $\overline{\Het(P,f)}=\Hom(P,f)$.\label{i.1}
\item If two basic sets $K$ and $L$ are heteroclinically related, they are robustly so.
\end{enumerate}
\end{proposition}

We say that a heteroclinic class $\Het(P,f)$ is {\em mechanically without domination of index $i$} if it contains a mechanical obstruction to domination of index $i$. Note that we do not need now to consider $\Het(P,f)$ to be a pointed set, for this to be well-defined. 
 
 \cref{t.comprehensible}  then restates in terms of heteroclinic classes as follows:

\medskip

\noindent {\em There is a $C^r$-residual set of diffeomorphisms $f$, such that if a heteroclinic class  $\Het(P,f)$ is mechanically without domination of any index, then its closure is the homoclinic class $\Hom(P,f)$ and is the Hausdorff limit of a sequence of sinks and sources.}
\medskip

%
%
%

\section{The $C^{1+\alpha}$-topologies}\label{s.C1alpha}
Newhouse, Palis-Viana and Turaev worked in $C^r$ topology, for $r\geq 2$. Our work also holds for the $C^{1+\alpha}$-topologies; it actually even holds in some way for the $C^1$-topology, provided we restrict ourselves to $C^{1+\alpha}$-diffeomorphisms with a fixed Hölder constant.

Let us give a few definitions.
Choose an atlas on $M$ and let $0<\alpha\leq 1$. Given a real number $C>0$, we denote by $\Diff^{1+\alpha}_C(M)$ the set of $C^1$-diffeomorphisms $f$ whose derivative has class $C^\alpha$ with Hölder constant $C$: for all points $x,y\in M$ such that $x,y$ are in a common chart and $f(x),f(y)$ also, we have
$$\|Df(x)-Df(y)\|<C\|x-y\|^\alpha,$$
for the respective norms given by these charts.
We endow that space with the topology of  $C^1$-convergence. This is a Baire space.

One of the many ingredients of \cref{t.robustization} is actually to turn robust homoclinic tangencies. This was done in $C^r$-topologies for $r\geq 2$, by Romero \cite{Ro} using the result of Palis and Viana~\cite{PV}. 

A result of Moreira~\cite{Mo} implies that this is not possible in the $C^1$-topology.
However, using a theorem of Bonatti and Crovisier \cite{BC2}, and adapting the original argument of Newhouse in dimension $2$, Crovisier and the author~\cite{CG} simplify and extend the Palis-Viana-Romero result about robustization of tangencies  to the space of $C^1$-diffeomorphisms with bounded Hölder constant.\footnote{This means that Moreira, when disjoining dynamically defined Cantor sets, builds $C^1$-maps that have huge Hölder constants if we coerce them to the $C^{1+\alpha}$-world.}
One then gets the following $C^{1+\alpha}$-version of  \cref{t.robustization}:
%
%

\begin{theorem}\label{t.robustizationalpha}
Fix an integer $r\geq 1$ and a real number $0<\alpha\leq 1$. Let $\Lambda$ be a homoclinic tangency or a heteroclinic cycle containing a saddle $P$ for a diffeomorphism $f\in \Diff^r(M)$. Assume that $\Lambda$ has no dominated splitting of any index in a set $I\subset \NN$. Let $C>0$ and let $\cV\subset \Diff^r(M)$ be a neighborhood of $f$.

There exists an open set $\cU\subset \Diff^1(M)$ such that
\begin{itemize}
\item $\cU$ intersects $\cV$,
\item if $g\in\cU\cap \Diff^{1+\alpha}_C(M)$, that is, if $g\in \cU$ is $C^{1+\alpha}$ with Hölder constant $C$, then the continuated chain-recurrent class $\cC(P_g,g)$ of $P$ has mechanically no domination of any index $i\in I$.
\end{itemize}
\end{theorem}

Then, using an adapted version of \cref{t.comprehensible}, one obtains the following $C^{1+\alpha}$-version of \cref{t.lambda}:

\begin{theorem}\label{t.lambdaalpha}
Fix real numbers $C>0$ and $0<\alpha<1$. Let $\Lambda$ be a homoclinic tangency or a heteroclinic cycle for a diffeomorphism $f\in \Diff^{1+\alpha}_C(M)$. One has a sharp dichotomy:
\begin{itemize}
\item either $\Lambda$ is volume-hyperbolic,
\item or there is an open set $\cU\subset \Diff^{1+\alpha}_C(M)$ such that 
 \begin{itemize}
 \item $f$ is in the closure of $\cU$
 \item and any $g$ in a residual set $\cR\subset \cU$ has infinitely many sinks or sources.
\end{itemize}
\end{itemize}
\end{theorem}

Actually, the set of diffeomorphisms in $\cU$ that have infinitely many sinks or sources is slightly larger than $C^1$-residual: it contains a countable intersection of $C^1$-open and, for all $r>1$, $C^r$-dense sets.

\section{Mechanical non-volume-hyperbolicity and Newhouse phenomena}\label{s.volhyp}

The most general result of~\cite{BDP} states that for any diffeomorphism $f$ in a residual subset of $\Diff^1(M)$, if a homoclinic class $\Hom(P,f)$ is not volume-hyperbolic, then it is the Hausdorff limit of a sequence of sinks or sources. This is in fact the theorem that we need to adapt to the $C^r$-topology in order to prove \cref{t.lambda}. Let us first define a notion of mechanical non-volume-hyperbolicity.

\begin{remark}
A compact invariant set $\Lambda$ of a diffeomorphism $f$ is not volume-hyperbolic if and only if, given the finest dominated splitting 

$$TM_{|\Lambda}=E^1\oplus ... \oplus E^k$$
above $\Lambda$, there exists a Riemannian metric on $M$ and a point $x\in \Lambda$ such that 
\begin{itemize}
\item either volume is not contracted on $E^1_x$ by any positive iterate of the dynamics, that is,  for any $n\in \NN$,  
 $$\bigl|\det Df^n_{|E^1_x}\bigr|\geq 1.$$
 \item or volume is not contracted on $E^k_x$ by any negative iterate of the dynamics, that is,  for any $n\in \NN$,  
  $$\bigl|\det Df^{-n}_{|E^k_x}\bigr|\geq 1.$$
\end{itemize}
\end{remark}

This remark motivates the following definition:

\begin{definition}
The chain-recurrent class $\Lambda=\cC(P,f)$ of a saddle $P$ for a diffeomorphism $f$ is {\em mechanically not volume-hyperbolic} if, given the finest dominated splitting 
\begin{align*}T_\Lambda M=E^1\oplus ... \oplus E^k\end{align*}
 above $\Lambda$, there exists a point $x$ in the heteroclinic class $\Het(P,f)$ of $P$ such that 
\begin{itemize}
\item either the bundle $E^1$ is mechanically without dominated splitting and the sequence  
$$\bigl|\det Df^n_{|E^1_x}\bigr|, \quad n\in \NN$$
is not bounded from above,
\item or the bundle $E^k$ is mechanically without dominated splitting and the sequence 
$$\bigl|\det Df^{-n}_{|E^k_x}\bigr|,  \quad n\in \NN$$
is not bounded from above.
\end{itemize}
\end{definition}
By the bundle $E^\ell$ being mechanically without dominated splitting, we mean that $\cC(P,f)$ has mechanically no dominations of all the corresponding indices, that is, of all indices $i$ such that 
$$\sum_{j<\ell} \dim E^j<i<\sum_{j\leq \ell} \dim E^j.$$
We conjecture that non-volume-hyperbolicity is generically mechanical:

\begin{conjecture}[]\label{c.volhyp}
For any integer $r\geq 1$, there is a residual subset $\cR$ of $C^r$-diffeomorphisms $f$ such that if the chain-recurrent class $\cC(f,P)$ of a saddle point $P$ is not volume-hyperbolic, then it is $C^r$-robustly mechanically not volume hyperbolic. 
\end{conjecture}

We now state the announced $C^r$-version of the most general theorem of~\cite{BDP}. 

\begin{theorem}\label{t.mechnovolhyp}
For any integer $r\geq 1$, for any diffeomorphism $f$ in a residual subset $\cR \subset \Diff^r(M)$, the following holds:

if a chain-recurrent class $\cC(P,f)$ is mechanically not volume-hyperbolic, then $\Hom(P,f)$ is the Hausdorff limit of a sequence of sinks or sources.
\end{theorem}

This straightforwardly implies~\cref{t.comprehensible}. Conjecture~\ref{c.volhyp} would then imply the following:

\begin{conjecture}
For any integer $r\geq 1$, for any diffeomorphism $f$ in a residual subset $\cR \subset \Diff^r(M)$ the following holds:

a homoclinic class $\Hom(P,f)$ is not volume-hyperbolic if and only if it is the Hausdorff limit of a sequence of sinks or sources.
\end{conjecture}

\section{Extensions of \cref{t.lambda,t.robustization} to other sets $\Lambda$}\label{s.moregenerallambda}

\cref{t.lambda,t.robustization}  solve \cref{p.newh} completely when the set $\Lambda$ is a tangency or a heteroclinic cycle. 
We explore in this section the problem for a wider category of sets $\Lambda$, namely unions of heteroclinically connected basic pieces.  We bring in particular a negative answer to the following:

\begin{question}\label{q.lambda}
Does the absence of volume-hyperbolicity on a compact chain-recurrent set $\Lambda$ imply that Newhouse phenomena occur $C^r$-close to $\Lambda$ ? 
\end{question}

The answer to this question, the extensions of \cref{t.lambda,t.robustization} and their addendum are not the only excuses for our introduction of the notion of heteroclinic graphs, later in this section; heteroclinic graphs also happen to be a convenient tool for understanding the proof of our $C^r$-versions of~\cite{BDP} and~\cite{BB}, they will also allow to quantify how complex the orbits of periodic points need to be, when one wants to prescribe such-and-such Lyapunov exponents, and later on, to prescribe such-and-such $n$-th order derivatives.

The next section answers negatively to \cref{q.lambda} in dimension $\geq 3$.

\subsection{An example of a non-dominated chain-recurrent set, far from Newhouse phenomena}\label{s.counterex}
We only build the counterexample in dimension $3$ as it is clear how to adapt it to greater dimension. Let $f$ be a diffeomorphism with two saddles $P$ and $Q$ of stable indices $1$ and $2$, respectively, such that the unstable manifold of $Q$ intersects that of $P$ at a point $x$, and the unstable manifold of $P$ intersects the stable one of $Q$ at two points $y_1$ and $y_2$. 

\begin{figure}
\begin{center}
\def\svgwidth{200pt}           
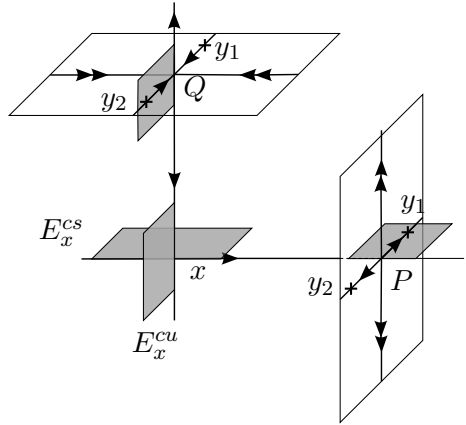
\end{center}
\caption{\small The chain-recurrent set $\Lambda=\{P,Q\}\cup \Orb(x)\cup \Orb(y_1)\cup \Orb(y_2)$ has no dominated splitting and is far from Newhouse phenomena. For the sake of clarity we chose not to represent the bundle tangencies at $y_1$ and $y_2$: both points appear twice in the figure.}
\label{f.counterex2}
\end{figure}

Assume moreover that 
\begin{itemize}
\item there is a three-bundles dominated splitting on $\{P,Q\}\cup \Orb(x)$,
\item there is a bundle tangency of index $i$, and only $i$, at the heteroclinic point $y_i$, for $i=1,2$,
\item the saddles $P$ and $Q$ are not sectionally dissipative, that is, if $\lambda_1<\lambda_2<\lambda_3$ are the Lyapunov exponents at $P$ (resp. $Q$), then
\begin{align*}
\lambda_1+\lambda_2&<0\\
\lambda_2+\lambda_3&>0.
\end{align*}
\item the points $y_i$ are separated by the strong-stable (resp. strong-unstable) manifold, within the stable (resp. unstable) manifold of $Q$ (resp. $P$).
\item if we denote by $E^{cs}_x$ the center-stable half plane in $T_xM$ whose positive iterates converge to the center-stable half plane at $P$ tangent to the half unstable manifold containing $y_1$, and by  $E^{cu}_x$ the center-unstable half plane at $x$ whose negative iterates converge to the center-unstable half plane at $Q$ tangent to the half unstable manifold containing $y_2$, then we have
$$E^{cs}_x\cap E^{cu}_x=\{x\}.$$
\end{itemize}
The example is depicted in  \cref{f.counterex2}.

\begin{claim}
There exist a neighborhood $\cU$ of $f$ in $\Diff^1(M)$ and a neigborhood $U$ in $M$ of the chain-recurrent set
$$\Lambda=\{P,Q\}\cup \Orb(x)\cup \Orb(y_1)\cup \Orb(y_2),$$ such that for any diffeomorphism $g\in \cU$, each periodic point in $U$ is neither a sink nor a source. 
\end{claim}

In particular, there is no Newhouse phenomena $C^r$-close to $\Lambda$.
The idea of the proof of the claim is as follows: if a periodic orbit $R$ is created close to $\Lambda$ by $C^1$-perturbations, then because of the geometry of $\Lambda$ it cannot pass both close to $y_1$ and close to $y_2$. As a consequence there is a strong domination of index either $1$ or $2$ along $R$. The absence of sectional dissipativity then prevents that $R$ be a sink or a source.

\subsection{\cref{t.lambda} for heteroclinic graphs}\label{s.newph}

We first need a few definitions. Recall that a {\em path} in a directed graph is an alternating sequence of vertices and connecting edges. A directed graph is called {\em strongly connected} if there exists a path from any vertex to any other.

\begin{definition}
A {\em heteroclinic graph} $\Gamma$ for a diffeomorphism $f$ is a directed graph, not necessarily finite, such that:
\begin{itemize}
\item each of its vertices is a chain-recurrent hyperbolic set,
\item each edge from vertex $K$ to vertex $L$ is an orbit whose $\alpha$-limit and $\omega$-limit are subsets of $K$ and $L$, respectively,
\end{itemize}
\end{definition}

A heteroclinic graph may contain loops (homoclinic orbits).
The {\em support} of a heteroclinic graph $\Gamma$, denoted by $\supp \Gamma$, is the union of the edges and vertices of $\Gamma$.

\begin{remark}
We allow two vertices to have nonempty intersection. This may bring some technical difficulties (see \cref{p.gathering}, for instance), but it allows the following nice example:
\end{remark}

\begin{example}
By definition, any heteroclinic class is the support of a strongly connected heteroclinic graph.
Note that the vertices of that heteroclinic graph cannot be chosen to be pairwise disjoint, in general.
\end{example}

A heteroclinic graph $\Gamma$ with support $\Lambda$ is {\em mechanically without domination of index $i$}  if so is 
$\Lambda$, that is, if $\Lambda$ contains a mechanical obstruction to domination of index $i$.
As seen in \cref{s.volhyp}, it is natural to say that $\Gamma$ is {\em mechanically not volume-hyperbolic} if, given the finest dominated splitting 
\begin{align*}T_{\Lambda} M=E^1\oplus ... \oplus E^k\end{align*}
there exists a point $x\in \Lambda$ such that 
\begin{itemize}
\item either the bundle $E^1$ is mechanically without domination and the sequence  
$$\bigl|\det Df^n_{|E^1_x}\bigr|, \quad n\in \NN$$
is not bounded from above,
\item or the bundle $E^k$ is mechanically without domination and the sequence 
$$\bigl|\det Df^{-n}_{|E^k_x}\bigr|,  \quad n\in \NN$$
is not bounded from above.
\end{itemize}

A {\em circuit} in a directed graph is a path that starts and ends at the same vertex, and that does not go twice through the same edge.
A directed graph is called {\em Eulerian} if for any finite subset of edges, there exists a circuit that goes through each of those edges. This implies in particular that the graph is strongly connected.

\begin{theorem}\label{t.volhyp}
 Let $r\geq 2$ be an integer. Let $\Lambda\subset M$ be the support of a Eulerian heteroclinic graph that is mechanically not volume-hyperbolic, then Newhouse phenomena occur $C^r$-close to $\Lambda$.
\end{theorem}

We have the following corollary:

\begin{corollary}\label{c.homcor}
 Let $r\geq 2$ be an integer and let $\Lambda\subset M$ be the support of a Eulerian heteroclinic graph that has mechanically no dominated splitting, then Newhouse phenomena occur $C^r$-close to $\Lambda$.
\end{corollary}

\begin{figure}
\begin{center}
\def\svgwidth{200pt}           
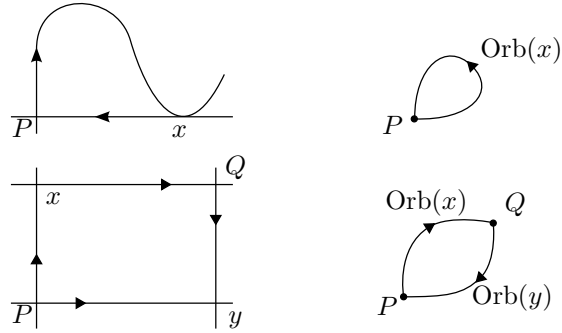
\end{center}
\caption{\small Homoclinic tangency and heterodimensional cycle with respective heteroclinic graphs. Both are Eulerian graphs.}
\label{f.graphecycletangence}
\end{figure}

The other all-important particular case is that if the vertices of a finite Eulerian heteroclinic graph are periodic points, then one has a sharp dichotomy between volume hyperbolicity and Newhouse phenomena: this is the addendum to \cref{t.lambda} in~\cref{s.turaevconj}. A homoclinic tangency and a heterodimensional cycle, as defined in \cref{s.turaevconj}, are the supports of finite Eulerian heteroclinic graphs (see \cref{f.graphecycletangence}), therefore we get \cref{t.lambda}.

The Eulerian hypothesis is actually generic:

\begin{proposition}\label{p.geneulerian}
There exists a residual set $\cR\subset \Diff^r(M)$ of diffeomorphisms $f$ such that any heteroclinic class $\Het(P,f)$ is the support of a Eulerian heteroclinic graph. 
\end{proposition}

\cref{p.geneulerian} and \cref{t.volhyp} allow to recover \cref{t.comprehensible,t.mechnovolhyp}.

\begin{figure}
\begin{center}
\def\svgwidth{100pt}           
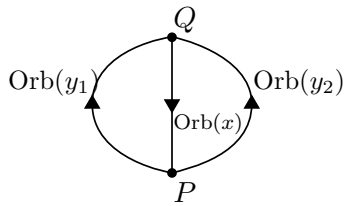
\end{center}
\caption{\small Strong connectivity of the graph is not sufficient: the support of this strongly connected but non-Eulerian graph is the counterexample of  \cref{s.counterex}.  It has mechanical obstructions to domination at all indices but no sinks or sources appear close to its support by $C^r$-perturbations.}
\label{f.counterex3}
\end{figure}

\begin{remark}One cannot replace the Eulerian hypothesis by the weaker strong connectivity hypothesis in  \cref{t.volhyp}, in dimension $\geq 3$. 
Indeed, the set $\Lambda=\{P,Q\}\cup \Orb(x)\cup \Orb(y_1)\cup \Orb(y_2)$ from the counterexample of \cref{s.counterex} is the support of a non-Eulerian strongly connected heteroclinic graph, depicted in \cref{f.counterex3}.
\end{remark}

\subsection{\cref{t.robustization}  for heteroclinic graphs}\label{s.robustization}

This part is about turning robust  mechanical lack of domination and is central to our work.

\subsubsection{Statement of result} A directed graph is {\em infinitely connected} if it is strongly connected and if it remains strongly connected by removal of any finite set of edges. We note that an infinitely connected graph is always Eulerian. We now state our result, which  straightforwardly implies \cref{t.robustization}:

\begin{theorem}\label{t.homcor}
Let $r\geq 2$ be an integer. If a saddle $P$ of $f\in \Diff^r(M)$ is in the support $\Lambda$ of a Eulerian heteroclinic graph that has mechanically no dominated splitting of indices $i\in I$, then for any open set $U$ containing $\Lambda$, there exists an open set $\cU\subset \Diff^r(M)$ such that:
\begin{itemize}
\item $f$ is in the closure of $\cU$,
\item for any $g$ in $\cU$, the continuation $P_g$ of the saddle $P$ is in an infinitely connected heteroclinic graph $\Delta_g$ supported in $U$, that contains infinitely many mechanical obstructions to domination of each index $i\in I$.
\end{itemize}
\end{theorem}

This is the central and most difficult part of our work. 

\begin{remark} The same way we wrote in \cref{s.C1alpha} a $C^{1+\alpha}$ version of \cref{t.robustization}, \cref{t.homcor} also admits a $C^{1+\alpha}$ version. We leave it to the reader to write its statement.
\end{remark}

\begin{remark} If the heteroclinic graph is finite and its vertices are saddles, one can remove the "mechanical" assumption in the hypotheses.
\end{remark}

\subsubsection{Ingredients for the proof of \cref{t.homcor}}
It is based on the next two propositions, which describe how heteroclinic graphs and obstructions to domination can be changed by $C^r$-perturbations.
A {\em trail} in a directed graph $\Gamma$ is a path that does not go twice through the same edge. A circuit is therefore a closed trail. The following propositions describe how Heteroclinic graphs and the mechanical obstructions at their vertices and edges may be transformed by $C^r$-perturbations of the dynamics. This will give us two "allowed" transformations on heteroclinic graphs. 

It would be of a clear interest to describe in more generality what are those allowed transformations. It is not even known when a cycle relating to saddles $P$ and $Q$ can be $C^r$-perturbed into a robust cycle containing both saddles (this was extensively studied by Bonatti and Diaz in the $C^1$-topology).

\begin{figure}
\begin{center}
\def\svgwidth{350pt}           
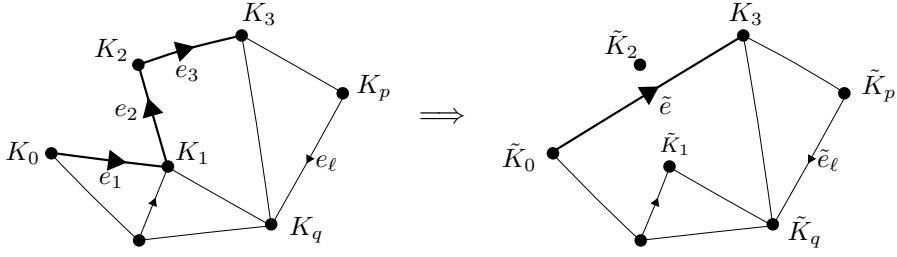
\end{center}
\caption{\small Reduction of a heteroclinic graph: all mechanical obstructions on $(e_1,K_1,e_2,K_2,e_3)$ are gathered on a single edge $\tilde{e}$ by a $C^r$-perturbation on a small neighborhood of $K_1\cup K_2$. The other vertices and edges retain their mechanical obstructions to domination of other indices.
\vskip 1mm
The reason why we need the Eulerian hypothesis in our theorems is that it is not possible, in general, to create the edge $\tilde{e}$ without destroying the edges $e_1,e_2,e_3$.}
\label{f.reductiongraphe}
\end{figure}

The first transformation that we present asserts that if a heteroclinic graph contains a trail $\gamma$ along which there are mechanical obstructions to domination of some indices $i\in I$, then by a $C^r$-perturbation one can gather those obstructions on a single edge with same origin and ending as the trail $\gamma$. 
The rest of the graph will keep the same configuration and mechanical obstructions for the indices $i\neq I$. Moreover, the support of the final graph can be found arbitrarily close to the support of the initial one. This transformation is illustrated in \cref{f.reductiongraphe}.

\begin{proposition}\label{p.gathering}
Let $r\geq 2$ be an integer. Let $f\in \Diff^r(M)$ be a diffeomorphism admitting a finite heteroclinic graph
$$\Gamma=\Bigl(\{K_0,...,K_m\},\{e_1,...,e_n\}\Bigr).$$
and assume that $\gamma=(K_0,e_1,K_1,\ldots K_{k-1},e_k,K_k)$ is a trail that contains mechanical obstructions to domination of any index $i$ in a set $I$. Let $\epsilon>0$ and $U$ a neighborhood of $K_1\cup \ldots \cup K_{k-1}$.

Then there is an $\epsilon$-perturbation $g$ of $f$ in $\Diff^r(M)$ and a heteroclinic graph for $g$
$$\tilde{\Gamma}=\Bigl(\{\tilde{K}_0,...,\tilde{K}_m\},\{\tilde{e},\tilde{e}_{k+1},...,\tilde{e}_n\}\Bigr)$$
such that $g=f$ outside $U$, and
\begin{itemize}
\item $(\tilde{K}_0,\tilde{e},\tilde{K}_k)$ is a trail that contains mechanical obstructions of all indices $i\in I$.
\item for all $0\leq \ell\leq m$, the basic set $\tilde{K}_\ell$ for $g$ is the continuation of the basic set $K_\ell$ for $f$ and retains its initial mechanical obstructions to domination, for all indices $i\notin I$.\footnote{Note that the vertices are not necessarily disjoint basic sets, which is why all the vertices, and not only vertices $K_1, \ldots ,K_{k-1}$, may loose domination of indices $i\in I$.}
\item for $\ell>k$, if $e_\ell$ is an edge from $K_p$ to $K_q$, then $\tilde{e}_\ell$ is an edge from $\tilde{K}_p$ to $\tilde{K}_q$ whose hausdorff distance to $e_\ell$ is less than $\epsilon$. Moreover, for any bundle tangency at $e_\ell$, there is a bundle tangency of same indices at $\tilde{e}_\ell$.
\end{itemize} 
\end{proposition}
In particular, the transformation can be done so that the Hausdorff distance between the supports of $\Gamma$ and  $\tilde{\Gamma}$ is arbitrarily small.

\begin{figure}
\begin{center}
\def\svgwidth{250pt}           
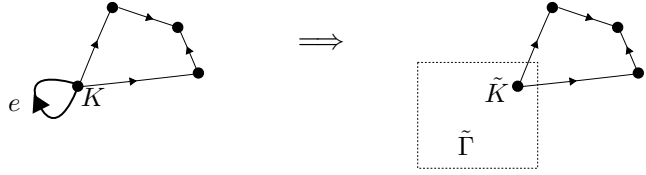
\end{center}
\caption{\small For any diffeomorphism $g$ in an open set close to the initial dynamics, the loop $(K,e,K)$ is becomes an infinitely connected heteroclinic graph $\tilde{\Gamma}$ containing the continuation $\tilde{K}$ of vertex $K$, and any mechanical obstruction in $(K,e,K)$ appears infinitely many times in $\tilde{\Gamma}$.
\vskip 1mm
If $K$ is part of another finite graph $\Delta$, then for one of those $g$, the graph $\Delta$ and its mechanical obstructions of other indices can survive in a continuated form.}
\label{f.completiongraphe}
\end{figure}

The second transformation allows to turn a loop  to a robustly infinitely connected heteroclinic graph so that any mechanical obstruction along the initial loop is infinitely replicated in the new graph.  The rest of the graph will keep the same configuration and mechanical obstructions of other indices on each vertex and edge (see \cref{f.completiongraphe}). The support of the final graph can be found arbitrarily close to the support of the initial one.

\begin{proposition}\label{p.robustizing}
Fix $r>1$ and $f\in \Diff^r(M)$. Let $\Gamma=\bigl(\{K\},\{e\}\bigr)$ be a heteroclinic graph for $f$ formed by a single loop $(K,e,K)$, and containing mechanical obstructions to domination of index $i$, for any integer $i$ in a set $I$.

Then there exists an open set $\cU\in \Diff^r(M)$ containing $f$ in its closure such that, for any $g\in \cU$,
\begin{itemize}
\item there is an infinitely connected heteroclinic graph $\tilde{\Gamma}$ for $g$ that contains the continuation $\tilde{K}$ of $K$,
\item  $\tilde{\Gamma}$ contains infinitely many mechanical obstructions to domination for each index $i\in I$,
\item the Hausdorff distance between the supports of $\Gamma$ and   $\tilde{\Gamma}$ is less than $\epsilon$.
\end{itemize} 
\end{proposition}

The conclusions of \cref{p.robustizing} may be complemented as follows:
\medskip

\noindent {\bf Addendum:}
{\em If $\Gamma$ is part of a larger finite heteroclinic graph $\Delta$, for any neighborhood $U$ of $K$ there is a diffeomorphism $g\in \cU$ arbitrarily $C^r$-close to $f$  such that 
\begin{itemize}
\item $g$ coincides with $f$ outside $U$,
\item $\Delta$ is still a heteroclinic graph for $g$, up to changing the vertices to their continuations, and up to adapting the edges adjacent to those continuated vertices, and the mechanical obstructions to domination of indices $i\notin I$ on each edge and each vertex are preserved.
\end{itemize} }

\cref{p.robustizing} is one of the most difficult result among those announced in this paper. Note that in the particular case of a homoclinic tangency in dimension $2$, this corresponds to the theorem of Newhouse that creates a robust homoclinic tangency by a $C^r$-perturbation.

\begin{proof}[Main steps of the proof of \cref{p.robustizing}:]
With \cref{p.gathering}, we reduce the problem to the case where $\Lambda$ is the closure of a homoclinic orbit $\Orb(x)$ related to a saddle $P$.  
\medskip

 Let $J$ be the set of indices of non-dominations that are visible only along the homoclinic orbit $\Orb(x)$, that is, through bundle tangencies at $x$. The  first step of the proof is to turn the homoclinic class of $P$ non-trivial, while recovering a homoclinic orbit $\Orb(y)$ of $P$ with bundle tangencies of same indices.
 
In dimension $2$, $J=\{1\}$, $\Lambda$ is a homoclinic tangency: a simple picture (see for instance~\cite{PT}) can convince the reader that a generic unfolding of the tangency leads to arbitrarily small parameters where the homoclinic class is non-trivial, and where the tangency is recovered. The same reasoning works in any dimension when $J=\{\ind P\}$: unfolding the tangency one recovers it for arbitrarily small parameters (see~\cite{PV,Ro}).
When the set $J$ has cardinal $\geq 2$, however, a generic unfolding will not allow to recover all those non-dominations at the same time. One needs to carefully choose the unfolding.

\medskip

Let now $I$ be the set of indices of non-domination along the closure of $\Orb(y)$. Each of those non-dominations is either due to a non-simple Lyapunov spectrum at $P$, that is, eigenvalues with same moduli at $P$, or a bundle tangency at $y$. By $C^r$-perturbations, one relates the bundle tangencies to other saddles homoclinically related to $P$. One is reduced to two independent problems:  given a non-trivial homoclinic class $H(P,f)$ of a saddle point $P$,
\begin{enumerate}
\item if $\lambda_1\leq ...\leq \lambda_d$ are the Lyapunov exponents along the saddle $P$ counted with multiplicity, find a $C^r$-perturbation of $f$ such that, for any index $i$, if $\lambda_i=\lambda_{i+1}$ then there is a saddle point $Q$ homoclinically related to the saddle $P$, at which there is robustly no dominated splitting of index $i$, thanks to a pair of complex conjugate eigenvalues.
\item given a bundle tangency related to $P$ of indices $i\in J$, find a $C^r$-perturbation such that
for any index $i\in J$, 
\begin{itemize}
\item if $i=\ind P$, there is a robust homoclinic tangency related to $P$.
\item if $i\neq \ind P$, there is a saddle point $Q$ homoclinically related to $P$ at which there is robustly no dominated splitting of index $i$, thanks to a pair of complex conjugate eigenvalues. 
\end{itemize}
\end{enumerate}
The first problem relies on the study of finitely generated semi-groups of matrices, and how those semi-groups behave through perturbations of their generators, according to the dominated splittings that exist on the corresponding linear cocycle.

The second problem is solved as follows: by perturbations of the eigenvalues at $P$ and perturbation close to the tangency, one manages to multiply bundle tangencies through  topological arguments. A bundle tangency of index $i=\ind P$, that is, a homoclinic tangency, can be turned robust thanks to~\cite{Ro} for regularity $r\geq 2$, or to~\cite{CG} for regularity $r>1$. Then a pair of bundle tangencies of same index $i\neq \ind(P)$ can be perturbed to obtain a saddle point homoclinically related to $P$, with complex conjugate eigenvalues at index $i$. 
\end{proof}

\section{Results \`a la Bochi-Bonatti in the $C^r$-topologies}\label{s.Lyapunov graphs}

In this section, we describe precisely the vectors of Lyapunov exponents one may obtain along periodic orbits by $C^r$-perturbations close to the chain-recurrent class of a point $P$, according to the mechanical obstructions to domination on it.
We actually give here a $C^r$-version of a theorem of $C^1$-dynamics by Bochi and Bonatti~\cite{BB} that generalizes~\cite{BDP}.

Let $\mu$ be an $f$-invariant probability measure with Lyapunov exponents $\lambda_1\leq \ldots \leq \lambda_d$, counted with multiplicity. The {\em Lyapunov map}\footnote{Bochi and Bonatti call it the {\em Lyapunov graph}, we save the word "graph" for heteroclinic graphs.} of $\mu$, identified to the tuple $\sigma(\mu)=(\sigma_0,\ldots,\sigma_d)$, where $\sigma_i=\sum_{1\leq j\leq i} \lambda_i$, is the map:
$$
\sigma(\mu)\colon \begin{cases}
\{0,\ldots,d\}&\to \mathbb{R}\\
i &\mapsto \sigma_i=\sum_{1\leq j\leq i} \lambda_j
\end{cases}.
$$
That map is convex by construction. For all $r\geq 1$ and $f\in \Diff^r(M)$, we say that an $f$-invariant probability measure $\mu$ is {\em $C^r$-preperiodic} if there exists a sequence $f_n\in \Diff^r(M)$ converging to $f$ and a sequence $P_n$ such that $P_n$ is a periodic orbit for $f_n$ and if the corresponding sequence of periodic measures $\mu_n$ converges weakly to $\mu$. 

Bochi and Bonatti get a complete description all the Lyapunov spectra one may obtain along periodic points by 
$C^1$-perturbations close to a $C^1$-preperiodic measure:
\medskip

\noindent {\bf Theorem} (Bochi-Bonatti~\cite[Theorem~3]{BB}){\bf .}
{\em Let $\mu$ be a $C^1$-preperiodic measure for $f$, with Lyapunov map $\sigma\colon \{0,\ldots,d\}\to \mathbb{R}$. Let $I$ be the set of integers $0\leq i\leq d$ such that there is a dominated splitting of index $i$ on the support of $\mu$ (we assume that $0,d\in I$). 

Then, any convex map $\tau\geq \sigma$ coinciding with $\sigma$ at indices $i\in I$ can be realized as a limit of Lyapunov maps for a sequence of diffeomorphisms $f_n$ converging to $f$ in $\Diff^1(M)$ and for a sequence of periodic measure $\mu_n\in \cM(f_n)$  weakly converging to  $\mu$.}
\medskip

\noindent The conclusion of the theorem may be more formally restated as follows: 
{\em For any convex map $\tau\geq \sigma$ coinciding with $\sigma$ at indices $i\in I$, there exists a sequence $f_n\in \Diff^1(M)$ converging to $f$ and a sequence $\mu_n$ of measures converging to $\mu$ for the weak-star topology such that:
\begin{itemize}
\item $\mu_n$ is a periodic measure for $f_n$ with Lyapunov map $\tau_n$, for all $n\in \NN$,
\item the sequence $\tau_n$ converges to $\tau$.
\end{itemize}}
\medskip

They actually have are more precise statement: they give a full description of such limit Lyapunov map, if they require moreover that the sequence of supports $\supp \mu_n$ converge for the Hausdorff topology to a fixed compact set $\Lambda$. One needs then to take for $I$ the set of integers $0\leq i\leq d$ such that there is a dominated splitting of index $i$ on $\Lambda$.

\begin{figure}
\begin{center}
\def\svgwidth{100pt}           
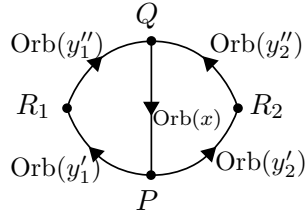
\end{center}
\caption{\small The sum of Dirac measures $\frac{1}{2}(\delta_{R_1}+\delta_{R_2})$ is a measure supported in a chain-recurrent set. However the periodic points that appear close to $\Lambda$, by perturbations, cannot visit both a small neighborhood of $R_1$ and of $R_2$.}
\label{f.counterex4}
\end{figure}

A central question for the global classification of $C^1$-generic dynamics is therefore to find when an invariant measure is $C^1$-preperiodic. Ma\~n\'e~\cite{M} proved that an ergodic measure is always $C^1$-preperiodic. 
Obviously, it is necessary that the measure be supported in a chain-recurrent set. This is however not sufficient as the following example shows: 

\begin{example}
The example of \cref{s.counterex} can be refined by adding a pair of saddles $R_1$ and $R_2$ along the heteroclinic connections $y_1$ and $y_2$: we get the heteroclinic graph $\Gamma$ depicted in \cref{f.counterex4}.  The sum of Dirac measures $\mu=\frac{1}{2}(\delta_{R_1}+\delta_{R_2})$ is a measure whose support is in the chain-recurrent set $\Lambda=\supp \Gamma$. 

However, for the same reasons as in  \cref{s.counterex}, a periodic point obtained by $C^1$-perturbations cannot go close to both $R_1$ and $R_2$. Thus $\mu$ is not $C^1$-preperiodic.
\end{example}

We conjecture however that:

\begin{conjecture}
For $r\geq 1$, there exists a $C^r$-residual set of diffeomorphisms for which any invariant measure supported in a chain-recurrent set is $C^r$-preperiodic. 
 \end{conjecture}

\subsection{Statement of results} We now give a statement similar to that of Bochi and Bonatti in the $C^r$ setting. The  invariant probability measure $\mu$ we now consider, instead of being $C^1$-preperiodic, is {\em supported in a heteroclinic class $\Het(P,f)$}, that is:
$$\mu\bigl[\Het(P,f)\bigr]=1.$$
Likewise, we say that a probability measure $\mu$ is {\em supported in a heteroclinic graph $\Gamma$} if the support of the graph has total measure:
 $$\mu(\supp \Gamma)=1.$$

\begin{theorem}\label{t.bobocr}
Let $\mu$ be an $f$-invariant probability measure supported in an infinitely connected heteroclinic graph $\Gamma$, with Lyapunov map $\sigma\colon \{0,\ldots,d\}\to \mathbb{R}$. Let $I$ be the set of integers $0\leq i\leq d$ such that there is no mechanical obstruction to domination of index $i$ on $\Gamma$ (we have $0,d\in I$). 

Then any convex map $\tau\geq \sigma$ coinciding with $\sigma$ at indices $i\in I$ can be realized as the limit Lyapunov map of a sequence of periodic measure $\mu_n$ converging weakly to $\mu$, for a sequence of diffeomorphisms $f_n$ converging to $f$ in $\Diff^r(M)$.

Moreover if $U$ is an open set containing the support of $\Gamma$, the sequence $(f_n,\mu_n)$ can be chosen so that $\supp \mu_n \subset U$ for $n$ large enough.
\end{theorem}

The infinite connection hypothesis is actually too strong. It would suffice that there are sufficiently many disjoint paths from a vertex to another; the needed number of such paths  only depends on the dimension $d$ of the manifold.
It is however an acceptable hypothesis since \cref{p.geneulerian} may be generalized to the following:

\begin{proposition}\label{p.geninfcon}
There exists a residual set $\cR\subset \Diff^r(M)$ of diffeomorphisms $f$ such that any heteroclinic class $\Het(P,f)$ is the support of an infinitely connected heteroclinic graph.
\end{proposition}

\cref{t.bobocr,p.geninfcon} then give the following generic consequence:

\begin{corollary}\label{c.Lyapgraph}
Let $r\geq 1$. There is a residual subset $\cR \subset \Diff^r(M)$ such that if 
\begin{itemize}
\item $f$ is a diffeomorphism in $\cR$,
\item $\mu$ is an $f$-invariant probability measure $\mu$ supported in a heteroclinic class $\Het(P,f)$, with Lyapunov map $\sigma\colon \{0,\ldots,d\}\to \mathbb{R}$,
\item $I$ is the set of integers $0\leq i\leq d$ such that there is no mechanical obstruction to domination of index $i$ on $\Het(P,f)$ -- we have $0,d\in I$ --
\end{itemize}
then any convex map $\tau\geq \sigma$ coinciding with $\sigma$ at indices $i\in I$ is the limit Lyapunov map of a sequence of periodic $f$-invariant measures $\mu_n$ converging weakly to $\mu$. 

Moreover, that sequence can be chosen so that the supports $\supp (\mu_n)$ converge to the homoclinic class $\Hom(P,f)$ for the Hausdorff topology.
\end{corollary}

This corollary straightforwardly implies \cref{t.comprehensible,t.mechnovolhyp}. 

\subsection{An application: polytope of realizable Lyapunov exponents along periodic orbits, close to a homoclinic tangency.}\label{s.polytope}

As is an example of an application of \cref{t.homcor,c.Lyapgraph}, we refine \cref{t.lambda} by giving a full description of the vectors of Lyapunov exponents of the periodic measures that can appear asymptotically close to a homoclinic tangency,  by $C^r$-perturbations:
these vectors form a convex polytope in an affine subspace of $\RR^d$. 

Remarkably, for $C^r$-generic diffeomorphisms $g$ in a nearby open set $\cU\subset \Diff^r(M)$, the vectors of Lyapunov exponents of the periodic measures of $g$ are dense in a (small) neighborhood of  that polytope.

\medskip

\begin{figure}
\begin{center}
\def\svgwidth{140pt}           
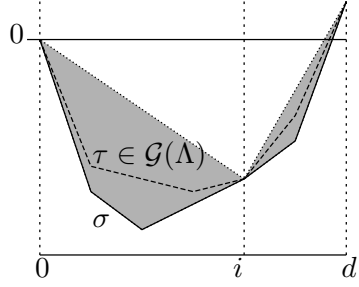
\end{center}
\caption{\small Assume that the only index of domination on $\Lambda$ is $i$, then $I=\{0,i,d\}$. The realizable limit Lyapunov maps close to $\Lambda$ are the convex maps in the gray region. They form a 3-dimensional polytope.}
\label{f.lyapunovmap}
\end{figure}

We first introduce some notation. Given a homoclinic tangency $\Lambda$ related to a saddle $P$ for a diffeomorphism $f$, let 
$$\sigma=\sigma(\delta_P)\colon \{0,\ldots,d\}\to \mathbb{R}$$ be the Lyapunov map 
for the Dirac measure at saddle point $P$, and let $I$ be the union of $\{0,d\}$ and the set of integers $0< i< d$ such that there is a dominated splitting of index $i$ on $\Lambda$. We denote by 
$$\cG(\Lambda)$$ 
the set of convex Lyapunov maps $\tau\geq \sigma$ coinciding with $\sigma$ at indices $i\in I$ (see~\cref{f.lyapunovmap}). The set $\cG(\Lambda)$ corresponds to a convex polytope inside an affine subspace of $\RR^d$ of codimension $\# (I)-1$.

\begin{theorem}\label{t.Lyaptang}
Let $r\geq 2$. Let $\Lambda$ be a homoclinic tangency for a diffeomorphism $f\in \Diff^r(M)$.
Then $\cG(\Lambda)$ is precisely the set of possible limit Lyapunov maps: $\tau\in \cG(\Lambda)$ if and only if there are
\begin{itemize}
\item a sequence of diffeomorphisms $f_n$ converging to $f$ in $\Diff^r(M)$
\item a sequence of periodic measures $\mu_n\in \cM(f_n)$ where the supports $\supp \mu_n$ converge to a subset of $\Lambda$ for the Hausdorff topology,
\end{itemize}
such that $\tau=\lim_{n\to \infty}\sigma(\mu_n)$. Moreover, for any open set  $U\supset \Lambda$, there exists an open set $\cU\subset \Diff^r(M)$
such that:
\begin{itemize}
\item $f$ is in the closure $\overline{\cU}$,
\item for any diffeomorphism $g$ in a residual subset $\cR\subset \cU$, if $\Pi_{g_{|U}}$ is the set of Lyapunov maps of periodic measures for the restriction $g_{|U}$, then the closure $\overline{\Pi_{g_{|U}}}$ is a neighborhood  in $\RR^d$ of the polytope $\cG(\Lambda)$.
\end{itemize}
\end{theorem}


\begin{remark}
\cref{t.Lyaptang} also holds
\begin{itemize}
\item if $\Lambda$ is the closure of the orbit of a homoclinic point for $P$,
\item in the $C^{1+\alpha}$ topologies. One needs however to adapt the statement the same way we adapted the statement of \cref{t.lambda} to obtain \cref{t.lambdaalpha}.
\end{itemize}
\end{remark}

If $\Lambda$ is now a heterodimensional cycle linking two saddles $P$ and $Q$, $I$ is the union of $\{0,d\}$ and the set of integers $0< i< d$ such that there is a dominated splitting of index $i$ on $\Lambda$, and $\mu,\sigma$ are the Lyapunov maps of the Dirac measures at the two saddles $P$ and $Q$ of the homoclinic cycle, then the candidate set $\cG(\Lambda)$ in \cref{t.Lyaptang} would be the polytope formed  of convex Lyapunov maps $\tau$ such that, for some $0\leq \lambda\leq 1$,
\begin{itemize}
\item $\tau\geq \tau_\lambda=\lambda\mu+(1-\lambda)\sigma$,
\item $\tau$ coincides with $\tau_\lambda$ at indices $i\in I$.
\end{itemize} 

\begin{conjecture}\label{t.Lyaphet}
\cref{t.Lyaptang} still works if one replaces the homoclinic tangency by a heterodimensional cycle.
\end{conjecture}

It would work when the difference of indices between the two saddles is $1$, by creating a robust heterodimensional cycle relating both saddles. It seems indeed reasonable to adapt to the $C^r$-topology the results of Bonatti-Diaz-Kiriki~\cite{BDK}. 

However this does not work for a larger difference between indices: there are examples of heterodimensional cycles that cannot be turned robust by $C^1$-perturbations. 
 One would need to find another way to create a non-trivial homoclinic class $\Hom(P,f)$ where there are periodic points that spend most of their time close to $Q$. We do not know if this is possible yet, even by $C^1$-perturbations.

\begin{ack}
This work owes a great deal to conversations with Enrique Pujals who first suggested in 2011 that one should try to generalize to the $C^r$-topology a number of $C^1$-results, bringing to the knowledge of the author the works of Gonchenko, Shil'nikov and Turaev. 

The author would like to thank Dmitry Turaev for his enthusiastic support and conversations, the Universit\'e d'Orsay where the idea of this work was born, and Sylvain Crovisier and Christian Bonatti who carefully listened and brought their commanding understanding of these problems.
\end{ack}


\end{document}

%% file: cyclebasicsets.pdf_tex
\begingroup%
  \makeatletter%
  \providecommand\color[2][]{%
    \errmessage{(Inkscape) Color is used for the text in Inkscape, but the package 'color.sty' is not loaded}%
    \renewcommand\color[2][]{}%
  }%
  \providecommand\transparent[1]{%
    \errmessage{(Inkscape) Transparency is used (non-zero) for the text in Inkscape, but the package 'transparent.sty' is not loaded}%
    \renewcommand\transparent[1]{}%
  }%
  \providecommand\rotatebox[2]{#2}%
  \ifx\svgwidth\undefined%
    \setlength{\unitlength}{271.97587891bp}%
    \ifx\svgscale\undefined%
      \relax%
    \else%
      \setlength{\unitlength}{\unitlength * \real{\svgscale}}%
    \fi%
  \else%
    \setlength{\unitlength}{\svgwidth}%
  \fi%
  \global\let\svgwidth\undefined%
  \global\let\svgscale\undefined%
  \makeatother%
  \begin{picture}(1,0.9376914)%
    \put(0,0){\includegraphics[width=\unitlength]{cyclebasicsets.pdf}}%
    \put(0.23570227,0.30392736){\color[rgb]{0,0,0}\makebox(0,0)[lb]{\smash{\small $P$}}}%
    \put(-0.0039066,0.4586179){\color[rgb]{0,0,0}\makebox(0,0)[lb]{\smash{$K_1$}}}%
    \put(0.85961535,0.01109932){\color[rgb]{0,0,0}\makebox(0,0)[lb]{\smash{$K_2$}}}%
    \put(0.68312911,0.90193455){\color[rgb]{0,0,0}\makebox(0,0)[lb]{\smash{$K_3$}}}%
    \put(0.0360129,0.0783321){\color[rgb]{0,0,0}\makebox(0,0)[lb]{\smash{$x_1$}}}%
    \put(0.85961535,0.49013332){\color[rgb]{0,0,0}\makebox(0,0)[lb]{\smash{$x_2$}}}%
    \put(0.27132789,0.8431058){\color[rgb]{0,0,0}\makebox(0,0)[lb]{\smash{$x_3$}}}%
  \end{picture}%
\endgroup%

%% file: grapheeulerienselles.pdf_tex
\begingroup%
  \makeatletter%
  \providecommand\color[2][]{%
    \errmessage{(Inkscape) Color is used for the text in Inkscape, but the package 'color.sty' is not loaded}%
    \renewcommand\color[2][]{}%
  }%
  \providecommand\transparent[1]{%
    \errmessage{(Inkscape) Transparency is used (non-zero) for the text in Inkscape, but the package 'transparent.sty' is not loaded}%
    \renewcommand\transparent[1]{}%
  }%
  \providecommand\rotatebox[2]{#2}%
  \ifx\svgwidth\undefined%
    \setlength{\unitlength}{520.8875bp}%
    \ifx\svgscale\undefined%
      \relax%
    \else%
      \setlength{\unitlength}{\unitlength * \real{\svgscale}}%
    \fi%
  \else%
    \setlength{\unitlength}{\svgwidth}%
  \fi%
  \global\let\svgwidth\undefined%
  \global\let\svgscale\undefined%
  \makeatother%
  \begin{picture}(1,0.49377082)%
    \put(0,0){\includegraphics[width=\unitlength]{grapheeulerienselles.pdf}}%
    \put(0.34916369,0.08782913){\color[rgb]{0,0,0}\makebox(0,0)[lb]{\smash{\small $P_3$}}}%
    \put(-0.00407958,0.21069635){\color[rgb]{0,0,0}\makebox(0,0)[lb]{\smash{\SMALL $\Orb(x_1)$}}}%
    \put(0.33380528,0.4564307){\color[rgb]{0,0,0}\makebox(0,0)[lb]{\smash{\small $P_1$}}}%
    \put(0.85599098,0.4564307){\color[rgb]{0,0,0}\makebox(0,0)[lb]{\smash{\small $P_2$}}}%
    \put(0.84063257,0.08782913){\color[rgb]{0,0,0}\makebox(0,0)[lb]{\smash{\small $P_4$}}}%
    \put(0.5529025,0.33788481){\color[rgb]{0,0,0}\makebox(0,0)[lb]{\smash{\SMALL $\Orb(x_1)$}}}%
  \end{picture}%
\endgroup%

%% file: contreexempledynamiquesimple.pdf_tex
\begingroup%
  \makeatletter%
  \providecommand\color[2][]{%
    \errmessage{(Inkscape) Color is used for the text in Inkscape, but the package 'color.sty' is not loaded}%
    \renewcommand\color[2][]{}%
  }%
  \providecommand\transparent[1]{%
    \errmessage{(Inkscape) Transparency is used (non-zero) for the text in Inkscape, but the package 'transparent.sty' is not loaded}%
    \renewcommand\transparent[1]{}%
  }%
  \providecommand\rotatebox[2]{#2}%
  \ifx\svgwidth\undefined%
    \setlength{\unitlength}{413.85390625bp}%
    \ifx\svgscale\undefined%
      \relax%
    \else%
      \setlength{\unitlength}{\unitlength * \real{\svgscale}}%
    \fi%
  \else%
    \setlength{\unitlength}{\svgwidth}%
  \fi%
  \global\let\svgwidth\undefined%
  \global\let\svgscale\undefined%
  \makeatother%
  \begin{picture}(1,0.86522833)%
    \put(0,0){\includegraphics[width=\unitlength]{contreexempledynamiquesimple.pdf}}%
    \put(0.56058429,0.2512964){\color[rgb]{0,0,0}\makebox(0,0)[lb]{\smash{$y_2$}}}%
    \put(0.73455873,0.40594034){\color[rgb]{0,0,0}\makebox(0,0)[lb]{\smash{$y_1$}}}%
    \put(0.17397444,0.59924527){\color[rgb]{0,0,0}\makebox(0,0)[lb]{\smash{$y_2$}}}%
    \put(0.38660986,0.69589774){\color[rgb]{0,0,0}\makebox(0,0)[lb]{\smash{$y_1$}}}%
    \put(0.32861838,0.61857577){\color[rgb]{0,0,0}\makebox(0,0)[lb]{\smash{$Q$}}}%
    \put(0.23196591,0.13531344){\color[rgb]{0,0,0}\makebox(0,0)[lb]{\smash{$E^{cu}_x$}}}%
    \put(0.05799148,0.34794887){\color[rgb]{0,0,0}\makebox(0,0)[lb]{\smash{$E^{cs}_x$}}}%
    \put(0.34063953,0.27198064){\color[rgb]{0,0,0}\makebox(0,0)[lb]{\smash{$x$}}}%
    \put(0.71522824,0.2512964){\color[rgb]{0,0,0}\makebox(0,0)[lb]{\smash{$P$}}}%
  \end{picture}%
\endgroup%

%% file: graphecycletangence.pdf_tex
\begingroup%
  \makeatletter%
  \providecommand\color[2][]{%
    \errmessage{(Inkscape) Color is used for the text in Inkscape, but the package 'color.sty' is not loaded}%
    \renewcommand\color[2][]{}%
  }%
  \providecommand\transparent[1]{%
    \errmessage{(Inkscape) Transparency is used (non-zero) for the text in Inkscape, but the package 'transparent.sty' is not loaded}%
    \renewcommand\transparent[1]{}%
  }%
  \providecommand\rotatebox[2]{#2}%
  \ifx\svgwidth\undefined%
    \setlength{\unitlength}{1001.2125bp}%
    \ifx\svgscale\undefined%
      \relax%
    \else%
      \setlength{\unitlength}{\unitlength * \real{\svgscale}}%
    \fi%
  \else%
    \setlength{\unitlength}{\svgwidth}%
  \fi%
  \global\let\svgwidth\undefined%
  \global\let\svgscale\undefined%
  \makeatother%
  \begin{picture}(1,0.61428379)%
    \put(0,0){\includegraphics[width=\unitlength]{graphecycletangence.pdf}}%
    \put(0.30090632,0.3617646){\color[rgb]{0,0,0}\makebox(0,0)[lb]{\smash{\small $x$}}}%
    \put(0.00163471,0.3617646){\color[rgb]{0,0,0}\makebox(0,0)[lb]{\smash{\small $P$}}}%
    \put(0.69512806,0.3675646){\color[rgb]{0,0,0}\makebox(0,0)[lb]{\smash{\small $P$}}}%
    \put(0.87893429,0.51139021){\color[rgb]{0,0,0}\makebox(0,0)[lb]{\smash{\small $\Orb(x)$}}}%
    \put(0.001640613,0.00985009){\color[rgb]{0,0,0}\makebox(0,0)[lb]{\smash{\small $P$}}}%
    \put(0.40099212,0.29235028){\color[rgb]{0,0,0}\makebox(0,0)[lb]{\smash{\small $Q$}}}%
    \put(0.06522044,0.23585023){\color[rgb]{0,0,0}\makebox(0,0)[lb]{\smash{\small $x$}}}%
    \put(0.40744923,0.01146434){\color[rgb]{0,0,0}\makebox(0,0)[lb]{\smash{\small $y$}}}%
    \put(0.6851456216,0.02512269){\color[rgb]{0,0,0}\makebox(0,0)[lb]{\smash{\small $P$}}}%
    \put(0.92459319,0.22145602){\color[rgb]{0,0,0}\makebox(0,0)[lb]{\smash{\small $Q$}}}%
    \put(0.701113775,0.22574868){\color[rgb]{0,0,0}\makebox(0,0)[lb]{\smash{\small $\Orb(x)$}}}%
    \put(0.86295367,0.04594244){\color[rgb]{0,0,0}\makebox(0,0)[lb]{\smash{\small $\Orb(y)$}}}%
  \end{picture}%
\endgroup%

%% file: counterexample.pdf_tex
\begingroup%
  \makeatletter%
  \providecommand\color[2][]{%
    \errmessage{(Inkscape) Color is used for the text in Inkscape, but the package 'color.sty' is not loaded}%
    \renewcommand\color[2][]{}%
  }%
  \providecommand\transparent[1]{%
    \errmessage{(Inkscape) Transparency is used (non-zero) for the text in Inkscape, but the package 'transparent.sty' is not loaded}%
    \renewcommand\transparent[1]{}%
  }%
  \providecommand\rotatebox[2]{#2}%
  \ifx\svgwidth\undefined%
    \setlength{\unitlength}{375.0875bp}%
    \ifx\svgscale\undefined%
      \relax%
    \else%
      \setlength{\unitlength}{\unitlength * \real{\svgscale}}%
    \fi%
  \else%
    \setlength{\unitlength}{\svgwidth}%
  \fi%
  \global\let\svgwidth\undefined%
  \global\let\svgscale\undefined%
  \makeatother%
  \begin{picture}(1,0.70069952)%
    \put(0,0){\includegraphics[width=\unitlength]{counterexample.pdf}}%
    \put(0.4905523,-0.02132834){\color[rgb]{0,0,0}\makebox(0,0)[lb]{\smash{$P$}}}%
    \put(0.4905522,0.63985074){\color[rgb]{0,0,0}\makebox(0,0)[lb]{\smash{$Q$}}}%
    \put(-0.12797014,0.40523882){\color[rgb]{0,0,0}\makebox(0,0)[lb]{\smash{\small $\Orb(y_1)$}}}%
    \put(0.7891492,0.40523882){\color[rgb]{0,0,0}\makebox(0,0)[lb]{\smash{\small$\Orb(y_2)$}}}%
    \put(0.4905522,0.25594032){\color[rgb]{0,0,0}\makebox(0,0)[lb]{\smash{\SMALL $\Orb(x)$}}}%
  \end{picture}%
\endgroup%

%% file: reductiongraphe.pdf_tex
\begingroup%
  \makeatletter%
  \providecommand\color[2][]{%
    \errmessage{(Inkscape) Color is used for the text in Inkscape, but the package 'color.sty' is not loaded}%
    \renewcommand\color[2][]{}%
  }%
  \providecommand\transparent[1]{%
    \errmessage{(Inkscape) Transparency is used (non-zero) for the text in Inkscape, but the package 'transparent.sty' is not loaded}%
    \renewcommand\transparent[1]{}%
  }%
  \providecommand\rotatebox[2]{#2}%
  \ifx\svgwidth\undefined%
    \setlength{\unitlength}{1013.84375bp}%
    \ifx\svgscale\undefined%
      \relax%
    \else%
      \setlength{\unitlength}{\unitlength * \real{\svgscale}}%
    \fi%
  \else%
    \setlength{\unitlength}{\svgwidth}%
  \fi%
  \global\let\svgwidth\undefined%
  \global\let\svgscale\undefined%
  \makeatother%
  \begin{picture}(1,0.25303234)%
    \put(0,0){\includegraphics[width=\unitlength]{reductiongraphe.pdf}}%
    \put(-0.00104799,0.09351562){\color[rgb]{0,0,0}\makebox(0,0)[lb]{\smash{\small $K_0$}}}%
    \put(0.18043954,0.09351562){\color[rgb]{0,0,0}\makebox(0,0)[lb]{\smash{\small $K_1$}}}%
    \put(0.2514564,0.2434401){\color[rgb]{0,0,0}\makebox(0,0)[lb]{\smash{\small $K_3$}}}%
    \put(0.09877015,0.0676222){\color[rgb]{0,0,0}\makebox(0,0)[lb]{\smash{\small $e_1$}}}%
    \put(0.09364116,0.20398629){\color[rgb]{0,0,0}\makebox(0,0)[lb]{\smash{\small $K_2$}}}%
    \put(0.115268,0.14086019){\color[rgb]{0,0,0}\makebox(0,0)[lb]{\smash{\small $e_2$}}}%
    \put(0.18043954,0.18320477){\color[rgb]{0,0,0}\makebox(0,0)[lb]{\smash{\small $e_3$}}}%
    \put(0.52974232,0.08562486){\color[rgb]{0,0,0}\makebox(0,0)[lb]{\smash{\small $\tilde{K}_0$}}}%
    \put(0.69886262,0.10219546){\color[rgb]{0,0,0}\makebox(0,0)[lb]{\smash{\SMALL $\tilde{K}_1$}}}%
    \put(0.77224671,0.2434401){\color[rgb]{0,0,0}\makebox(0,0)[lb]{\smash{\small $\tilde{K}_3$}}}%
    \put(0.69795966,0.143369){\color[rgb]{0,0,0}\makebox(0,0)[lb]{\smash{\small $\tilde{e}$}}}%
    \put(0.4408347,0.13296943){\color[rgb]{0,0,0}\makebox(0,0)[lb]{\smash{$\Longrightarrow$}}}%
    \put(0.37432685,0.16565971){\color[rgb]{0,0,0}\makebox(0,0)[lb]{\smash{\small $K_p$}}}%
    \put(0.30218273,0.01235343){\color[rgb]{0,0,0}\makebox(0,0)[lb]{\smash{\small $K_q$}}}%
    \put(0.33036403,0.08562486){\color[rgb]{0,0,0}\makebox(0,0)[lb]{\smash{\small $e_\ell$}}}%
    \put(0.86693586,0.08562486){\color[rgb]{0,0,0}\makebox(0,0)[lb]{\smash{\small $\tilde{e}_\ell$}}}%
    \put(0.83537281,0.00671714){\color[rgb]{0,0,0}\makebox(0,0)[lb]{\smash{\small $\tilde{K}_q$}}}%
    \put(0.91428043,0.16453248){\color[rgb]{0,0,0}\makebox(0,0)[lb]{\smash{\small $\tilde{K}_p$}}}%
    \put(0.63923102,0.20736806){\color[rgb]{0,0,0}\makebox(0,0)[lb]{\smash{\small $\tilde{K}_2$}}}%
  \end{picture}%
\endgroup%

%% file: completiongraphe.pdf_tex
\begingroup%
  \makeatletter%
  \providecommand\color[2][]{%
    \errmessage{(Inkscape) Color is used for the text in Inkscape, but the package 'color.sty' is not loaded}%
    \renewcommand\color[2][]{}%
  }%
  \providecommand\transparent[1]{%
    \errmessage{(Inkscape) Transparency is used (non-zero) for the text in Inkscape, but the package 'transparent.sty' is not loaded}%
    \renewcommand\transparent[1]{}%
  }%
  \providecommand\rotatebox[2]{#2}%
  \ifx\svgwidth\undefined%
    \setlength{\unitlength}{801.47231445bp}%
    \ifx\svgscale\undefined%
      \relax%
    \else%
      \setlength{\unitlength}{\unitlength * \real{\svgscale}}%
    \fi%
  \else%
    \setlength{\unitlength}{\svgwidth}%
  \fi%
  \global\let\svgwidth\undefined%
  \global\let\svgscale\undefined%
  \makeatother%
  \begin{picture}(1,0.25322417)%
    \put(0,0){\includegraphics[width=\unitlength]{completiongraphe.pdf}}%
    \put(0.43216223,0.18036984){\color[rgb]{0,0,0}\makebox(0,0)[lb]{\smash{$\Longrightarrow$}}}%
    \put(-0.00132569,0.08768333){\color[rgb]{0,0,0}\makebox(0,0)[lb]{\smash{\small $e$}}}%
    \put(0.71164787,0.10051681){\color[rgb]{0,0,0}\makebox(0,0)[lb]{\smash{\small $\tilde{K}$}}}%
    \put(0.67190505,0.02066377){\color[rgb]{0,0,0}\makebox(0,0)[lb]{\smash{\small $\tilde{\Gamma}$}}}%
    \put(0.10721979,0.09392875){\color[rgb]{0,0,0}\makebox(0,0)[lb]{\smash{\small $K$}}}%
  \end{picture}%
\endgroup%

%% file: counterexample2.pdf_tex
\begingroup%
  \makeatletter%
  \providecommand\color[2][]{%
    \errmessage{(Inkscape) Color is used for the text in Inkscape, but the package 'color.sty' is not loaded}%
    \renewcommand\color[2][]{}%
  }%
  \providecommand\transparent[1]{%
    \errmessage{(Inkscape) Transparency is used (non-zero) for the text in Inkscape, but the package 'transparent.sty' is not loaded}%
    \renewcommand\transparent[1]{}%
  }%
  \providecommand\rotatebox[2]{#2}%
  \ifx\svgwidth\undefined%
    \setlength{\unitlength}{375.0875bp}%
    \ifx\svgscale\undefined%
      \relax%
    \else%
      \setlength{\unitlength}{\unitlength * \real{\svgscale}}%
    \fi%
  \else%
    \setlength{\unitlength}{\svgwidth}%
  \fi%
  \global\let\svgwidth\undefined%
  \global\let\svgscale\undefined%
  \makeatother%
  \begin{picture}(1,0.70069952)%
    \put(0,0){\includegraphics[width=\unitlength]{counterexample2.pdf}}%
    \put(0.42656712,-0.04265668){\color[rgb]{0,0,0}\makebox(0,0)[lb]{\smash{$P$}}}%
    \put(0.426567,0.6611791){\color[rgb]{0,0,0}\makebox(0,0)[lb]{\smash{$Q$}}}%
    \put(-0.04265671,0.08531346){\color[rgb]{0,0,0}\makebox(0,0)[lb]{\smash{\small $\Orb(y'_1)$}}}%
    \put(0.72516413,0.10664182){\color[rgb]{0,0,0}\makebox(0,0)[lb]{\smash{\small $\Orb(y'_2)$}}}%
    \put(0.4905522,0.27726867){\color[rgb]{0,0,0}\makebox(0,0)[lb]{\smash{\SMALL $\Orb(x)$}}}%
    \put(-0.02132836,0.31992539){\color[rgb]{0,0,0}\makebox(0,0)[lb]{\smash{$R_1$}}}%
    \put(0.85313427,0.31992539){\color[rgb]{0,0,0}\makebox(0,0)[lb]{\smash{$R_2$}}}%
    \put(-0.04265671,0.55453731){\color[rgb]{0,0,0}\makebox(0,0)[lb]{\smash{\small $\Orb(y''_1)$}}}%
    \put(0.70383577,0.55453731){\color[rgb]{0,0,0}\makebox(0,0)[lb]{\smash{\small $\Orb(y''_2)$}}}%
  \end{picture}%
\endgroup%

%% file: lyapunovmap.pdf_tex
\begingroup%
  \makeatletter%
  \providecommand\color[2][]{%
    \errmessage{(Inkscape) Color is used for the text in Inkscape, but the package 'color.sty' is not loaded}%
    \renewcommand\color[2][]{}%
  }%
  \providecommand\transparent[1]{%
    \errmessage{(Inkscape) Transparency is used (non-zero) for the text in Inkscape, but the package 'transparent.sty' is not loaded}%
    \renewcommand\transparent[1]{}%
  }%
  \providecommand\rotatebox[2]{#2}%
  \ifx\svgwidth\undefined%
    \setlength{\unitlength}{234.29375bp}%
    \ifx\svgscale\undefined%
      \relax%
    \else%
      \setlength{\unitlength}{\unitlength * \real{\svgscale}}%
    \fi%
  \else%
    \setlength{\unitlength}{\svgwidth}%
  \fi%
  \global\let\svgwidth\undefined%
  \global\let\svgscale\undefined%
  \makeatother%
  \begin{picture}(1,0.76229968)%
    \put(0,0){\includegraphics[width=\unitlength]{lyapunovmap.pdf}}%
    \put(0.29023396,0.324379){\color[rgb]{0,0,0}\makebox(0,0)[lb]{\smash{$\tau \in \cG(\Lambda)$}}}%
    \put(0.1365807,0.01707247){\color[rgb]{0,0,0}\makebox(0,0)[lb]{\smash{$0$}}}%
    \put(0.66583084,0.01707247){\color[rgb]{0,0,0}\makebox(0,0)[lb]{\smash{$i$}}}%
    \put(0.95606479,0.01707247){\color[rgb]{0,0,0}\makebox(0,0)[lb]{\smash{$d$}}}%
    \put(0.29023396,0.15365315){\color[rgb]{0,0,0}\makebox(0,0)[lb]{\smash{$\sigma$}}}%
    \put(0.06829036,0.64875833){\color[rgb]{0,0,0}\makebox(0,0)[lb]{\smash{$0$}}}%
  \end{picture}%
\endgroup%